\documentclass[10pt]{amsart}                          
\usepackage{epsfig}            
\usepackage{amssymb,amscd,bbm}     
\usepackage{amsthm,amsgen,amsmath}
\usepackage[dvips]{xypic}       

\xyoption{curve}                
\xyoption{rotate}               

\newcommand{\Z}{{\mathbb Z}}  
\newcommand{\Q}{{\mathbb Q}}
       
\newcommand{\Lie}{\mathcal{L}\!{\mathit ie\/}}    
\newcommand{\Eil}{\mathcal{E}\!{\mathit il\/}}    
\newcommand{\Gr}{\mathcal{G}\!{\mathit r\/}}      
\newcommand{\Tr}{\mathcal{T}\!{\mathit r\/}}      
\newcommand{\Ass}{\mathcal{A}\!{\mathit s\/}}	  

\newcommand{\A}{{\mathcal A}}
\newcommand{\G}{{\mathcal G}}
\newcommand{\E}{{\mathcal E}}
\newcommand{\C}{{\mathcal C}}
\newcommand{\eL}{{\mathcal L}}
\newcommand{\la}{\langle}
\newcommand{\ra}{\rangle}

\newcommand{\dgc}{{\text{\scshape dgcc}}}
\newcommand{\dgl}{{\text{\scshape dgla}}}
\newcommand{\dga}{{\text{\scshape dgca}}}
\newcommand{\dge}{{\text{\scshape dglc}}}
\newcommand{\dgg}{{\text{\scshape dggc}}}

\newcommand{\dgt}{{\text{\scshape dgta}}}
\newcommand{\gc}{{\text{\scshape gcc}}}
\newcommand{\gl}{{\text{\scshape gla}}}
\newcommand{\ga}{{\text{\scshape gca}}}
\def\gg{{\text{\scshape ggc}}}

\newcommand{\dg}{{\text{\scshape dg}}}
\newcommand{\g}{{\text{\scshape g}}}

\newcommand{\freeA}{{\mathbb A}}
\newcommand{\freeT}{{\mathbb T}}
\newcommand{\vsfreeT}{\overline{{\mathbb T}}}
\newcommand{\cofreeC}{{\mathbb C}}
\newcommand{\cofreeE}{{\mathbb E}}
\newcommand{\freeL}{{\mathbb L}}
\newcommand{\cofreeG}{{\mathbb G}}
\newcommand{\vscofreeG}{{\overline{\mathbb G}}}
\newcommand{\vscofreeE}{{\overline{\mathbb E}}}
\newcommand{\vsfreeL}{{\overline{\mathbb L}}}
\newcommand{\gras}{\mathcal{A}\!{\mathit s\/}\mathcal{G}\!{\mathit r\/}}
\newcommand{\grac}{\mathcal{A}\mathcal{C}\mathcal{G}\!{\mathit r\/}}

\newcommand{\lie}{\mathcal{L}\!{\mathit ie\/}}    
\newcommand{\eil}{\mathcal{E}\!{\mathit il\/}}    

\newcommand{\inv}{ {{\text -}1}}      

\newcommand{\overtie}[2]{             
\begin{aligned} \displaystyle         %
\operatornamewithlimits{              %
 \begin{aligned} #1                   %
 \end{aligned} }_{#2}                 %
\end{aligned} }

\newcommand{\treel}[3]{ \ensuremath{ 
 \begin{xy}                          %
   (1,1); (2,2)**\dir{-},            
   (0,2); (2,0)**\dir{-};            
   (4,2)**\dir{-},                   
   (2,0); (2,-1)**\dir{-},           %
   (0,3.2)*{\scriptstyle #1},        %
   (2,3.2)*{\scriptstyle #2},        %
   (4,3.2)*{\scriptstyle #3},        %
 \end{xy}  } }

\newcommand{\longgraph}[4]{ \ensuremath{
 \begin{xy}                           
  (0,-2)*+UR{\scriptstyle #1}="a",    
  (4,3)*+UR{\scriptstyle #2}="b",     
  (8,-2)*+UR{\scriptstyle #3}="c",    
  (12,3)*+UR{\scriptstyle #4}="d",    %
  "a";"b"**\dir{-}?>*\dir{>},         %
  "b";"c"**\dir{-}?>*\dir{>},         %
  "c";"d"**\dir{-}?>*\dir{>}          %
 \end{xy}                             %
} }
\newcommand{\graphpp}[3]{ \ensuremath{
 \begin{xy}                           
  (0,-2)*+UR{\scriptstyle #1}="a",    
  (3,3)*+UR{\scriptstyle #2}="b",     
  (6,-2)*+UR{\scriptstyle #3}="c",    
  "a";"b"**\dir{-}?>*\dir{>},         %
  "b";"c"**\dir{-}?>*\dir{>}          %
 \end{xy}                             %
} }

\newcommand{\graphopp}[3]{ \ensuremath{%
 \begin{xy}                           
  (0,-2)*+UR{\scriptstyle #1}="a",    
  (3,3)*+UR{\scriptstyle #2}="b",     
  (6,-2)*+UR{\scriptstyle #3}="c",    
  "b";"c"**\dir{-}?>*\dir{>},         
  "c";"a"**\dir{-}?>*\dir{>}          %
 \end{xy}                             %
} }
\newcommand{\graphpop}[3]{ \ensuremath{%
 \begin{xy}                           
  (0,-2)*+UR{\scriptstyle #1}="a",    
  (3,3)*+UR{\scriptstyle #2}="b",     
  (6,-2)*+UR{\scriptstyle #3}="c",    
  "a";"b"**\dir{-}?>*\dir{>},         
  "c";"a"**\dir{-}?>*\dir{>}          %
 \end{xy}                             %
} }
\newcommand{\graphpm}[3]{ \ensuremath{
 \begin{xy}                           
  (0,-2)*+UR{\scriptstyle #1}="a",    
  (3,3)*+UR{\scriptstyle #2}="b",     
  (6,-2)*+UR{\scriptstyle #3}="c",    
  "a";"b"**\dir{-}?>*\dir{>},         %
  "c";"b"**\dir{-}?>*\dir{>}          %
 \end{xy}                             %
} }
\newcommand{\graphmp}[3]{ \ensuremath{
 \begin{xy}                           
  (0,-2)*+UR{\scriptstyle #1}="a",    
  (3,3)*+UR{\scriptstyle #2}="b",     
  (6,-2)*+UR{\scriptstyle #3}="c",    
  "b";"a"**\dir{-}?>*\dir{>},         
  "b";"c"**\dir{-}?>*\dir{>}          %
 \end{xy}                             %
} }

\newcommand{\linep}[2]{ \ensuremath{  %
 \begin{xy}                           
  (0,-2)*+UR{\scriptstyle #1}="a",    
  (3,3)*+UR{\scriptstyle #2}="b",     
  "a";"b"**\dir{-}?>*\dir{>},         
 \end{xy}                             %
} }
\newcommand{\smalllinep}[2]{ \ensuremath{  %
 \begin{xy}                           
  (0,-2)*+UR{\scriptstyle #1}="a",    
  (4,1)*+UR{\scriptstyle #2}="b",     
  "a";"b"**\dir{-}?>*\dir{>},         
 \end{xy}                             %
} }
\newcommand{\linem}[2]{ \ensuremath{  %
 \begin{xy}                           
  (0,-2)*+UR{\scriptstyle #1}="a",    
  (3,3)*+UR{\scriptstyle #2}="b",     
  "b";"a"**\dir{-}?>*\dir{>},         
 \end{xy}                             %
} }

\newcommand{\modcobr}[1]{             
 \raisebox{2pt}{$#1$}\!\diagup\,      %
 \raisebox{-2pt}{$\mathrm{ker}$ \raisebox{1pt}{$\scriptstyle ]\cdot[$}} }

\input colordvi                       

\theoremstyle{plain}                          
\newtheorem{theorem}{Theorem}[section]                          
\newtheorem{proposition}[theorem]{Proposition}                          
\newtheorem{lemma}[theorem]{Lemma}                          
\newtheorem{corollary}[theorem]{Corollary}

\theoremstyle{definition}                          
\newtheorem{definition}[theorem]{Definition}  

\theoremstyle{remark}                          
\newtheorem{remark}[theorem]{Remark}

\newtheorem{example}[theorem]{Example}

\newcommand{\refT}[1]{Theorem~\ref{T:#1}}
\newcommand{\refC}[1]{Corollary~\ref{C:#1}}
\newcommand{\refP}[1]{Proposition~\ref{P:#1}}
\newcommand{\refD}[1]{Definition~\ref{D:#1}}
\newcommand{\refL}[1]{Lemma~\ref{L:#1}}

\begin{document}

\title[Lie coalgebras I: graph coalgebras]
   {Lie coalgebras and rational homotopy theory, I: \\ graph coalgebras}

\author[D. Sinha]{Dev Sinha}
\address{Mathematics Department\\
University of Oregon\\
Eugene, OR
97403}
\email{dps@math.uoregon.edu}

\author[B. Walter]{Ben Walter} 
\address{
Department of Mathematics\\ Purdue University
150 N. University Street\\ West Lafayette, IN 47907
}
\email{walterb@math.purdue.edu}


\thanks{The first author was supported in part by NSF-DMS 0405922.\\
The second author was supported in part by the Mittag-Leffler Institute 
(Djursholm, Sweden),
and the Max-Planck Institute, Bonn.}

\date{\today}

\maketitle

\section{Introduction}

In this paper we develop a new, computationally friendly approach to Lie 
coalgebras through graph coalgebras,  and we apply this approach to 
Harrison homology.  There are two standard to presentations of a
Lie algebra through ``simpler'' algebras.  One is as a quotient of a
non-associative binary algebra by Jacobi and anti-commutativity identities.
Another presentation is as as embedded as Hopf algebra primitives in an
associative universal enveloping algebra.
The standard presentation of Lie coalgebras in the literature is dual to the second of these --
as a quotient of the associative coenveloping coalgebra, namely the 
Hopf algebra indecomposables \cite{Mich80, ScSt85}.
We describe an approach to Lie coalgebras
indiginous to the realm of coalgebras, dual to neither of these.
We define a new kind of coalgebra structure, namely anti-commutative
graph coalgebras, and we show that 
Lie coalgebras are quotients of these graph coalgebras.  

Our approach through graph coalgebras gives a presentation for Lie coalgebras 
which works better than the classical presentation 
in two respects.  First, cofree graph coalgebras come with a simple and easily computable 
pairing with free binary nonassociative
algebras which passes to Lie coalgebras and algebras,
 making duality not just a theoretical statement but
an explicitly computable tool.  Secondly, the quotient used to create Lie coalgebras
from graph coalgebras is a locally defined relation.  The quotient creating Lie
coalgebras from associative coalgebras is the shuffle relation, which causes global changes
to an expression.  As a result,  proofs in the realm of Lie 
coalgebras are often simpler to give through graph coalgebras than through associative
coalgebras, and for some important statements we have only found proofs in the graph coalgebra setting. 
For applications, 
we investigate the word problem for Lie coalgebras, and we also revisit 
Harrison homology.  
The category of graph coalgebras, and the graph cooperad on which it is based,
may also be of intrinsic interest.  The graph cooperad
is not binary, but could play a similar role in some  
natural category of cooperads as is played by the tree operad for binary operads.

\medskip

The plan of the paper is as follows.
After defining the graph cooperad, we pair it with the tree operad to give rise to a pairing 
on cofree and free algebras over them.  
We show that upon quotienting by the kernels of the pairing, it descends to 
a pairing between cofree Lie coalgebras and free Lie algebras.
This approach gives rise to our graphical model for the cofree Lie coalgebra
on a vector space $V$ and determines how that model pairs with the free Lie
algebra on a linear dual of $V$.  
Moreover, we can deduce a formula for the linear duality between Michaelis's 
Lie coalgebra model \cite{Mich80} and the tree/bracket model for free Lie algebras.  
We are also able to shed new light on the
structure of cofree Lie coalgebras, for example viewing them as what one
gets when one starts with a graph or associative coalgebra and ``kills the kernel
of the cobracket.''

We then lift the Andr\'e-Quillen construction on a 
differential graded commutative  algebra $(\dga)$
from the category of  differential graded Lie coalgebras $(\dge)$ to 
anti-commutative differential graded graph coalgebras
$(\dgg)$.  The Harrison model
for this bar construction passes through the category of associative coalgebras,
but our factorization through graph coalgebras is needed for example in 
developing an algebraic models for fibrations in the Lie
coalgebraic formulation of rational homotopy theory.   Such a result is 
critical in the sequel to this paper, where we define
generalized Hopf invariants and show from first principals that they give
a complete set of homotopy functionals in the simply-connected setting.
Indeed it was an investigation of
generalized Hopf invariants, which we found to be naturally indexed by
graphs, which led us to the framework of this paper.

Finally, we combine these results to shed new light on Quillen's
seminal work on rational homotopy theory \cite{Quil69}.  Quillen produced a pair of adjoint
functors $\eL$ and $\C$ between the categories of dg-commutative 
coalgebras $(\dgc)$ and dg-Lie algebras $(\dgl)$.  
In the linearly dual setting, there previously were two avenues 
towards understanding the functors between $\dga$ and $\dge$. 
One would be  a 
formal application of linear duality to Quillen's functors.   
The other way to  go from $\dga$ to $\dge$ explicitly was to use the Harrison complex,
which from Schlessinger and Stasheff's work has the structure of a Lie coalgebra dual
to Quillen's Lie algebraic functor.
Our techniques allow us to explicitly calculate the
linear duality between Harrison homology of a differential graded commutative algebra
and Quillen's functor $\eL$ on the corresponding linearly dual coalgebra, unifying these
approaches.

In our appendices, we take the opportunity to flesh out our models and connect with other work.
In particular, we give a spectral sequence for rational homotopy groups of a simply connected 
space,  we explicitly define model structures, and we discuss  minimal models.

\medskip

Our work throughout is over a field of characteristic zero.  We emphasize
that we are adding a finiteness hypothesis, namely that our algebras and
coalgebras are finite-dimensional in each positive degree, for the sake of
linear duality theorems.   Under this hypothesis the category
of chain complexes is canonically isomorphic to that of cochain complexes,
and we will use this isomorphism without further comment, by abuse denoting
both categories by $\dg$.  To clarify when possible, we have endeavored to 
use $V$ to denote a chain complex and $W$ to
denote a cochain complex.  Many of the facts we prove are true
without the finiteness hypothesis, as we may indicate.

We further restrict our work to 1-connected objects both to mirror the classical 
constructions of \cite{Quil69} and to allow ourselves to cleanly express our cofree Lie 
coalgebras as coinvariants rather than invariants. 
We plan to remove the finiteness and 
1-connectivity hypotheses in the third paper
in this series.
Note however that though in Sullivan's rational homotopy theory it is fairly typical to
quickly move to the nilpotent setting, this step requires a significant change to foundations
of our work.  The first author is currently writing a general theory of coalgebras over 
cooperads \cite{Walt08} so that we may proceed with such a program, where it looks like
we can extend even
beyond the nilpotent setting.

While we start by giving operadic definitions,  we work more
explicitly at the algebra and coalgebra level in later sections.  One reason for this
change in emphasis is a desire for explicit formulae.  But the change in emphasis is
necessary, since we have yet to find a purely operadic argument for the existence 
of the lift of the bar construction on a commutative algebra from the category of
Lie coalgebras to the category of graph coalgebras.
We hope to study the graph cooperad and graph coalgebras more extensively in future work.
We have yet to fully understand even what general (that is, not cofree)
graph coalgebras are in explicit algebraic terms.



\section{The graph cooperad and the configuration pairing}

We begin with constructions on the level of operads and cooperads,
to give more fundamental understanding (to readers familiar with operads) and
provide a general road-map for the following sections.
Later proofs and constructions will be given wholly in
the realm of algebras and coalgebras even when they could be inferred from operad
level statements presented here, which in some important cases they cannot be.  
A reader not interested in operads can skip most of this section, with the exceptions of
the definitions of graphs (\ref{D:graphs}), the 
configuration pairing (\ref{D:confpair}), and the quotients defining Lie coalgebras 
(\ref{D:eiln}).

\begin{definition}\label{D:graphs}  The graph symmetric sequence is defined as follows.
\begin{enumerate}
\item  Let $S$ be a finite set.
An $S$-graph is a connected oriented acyclic graph with vertex set 
$\mathrm{Vert}(G) = S$.  
\item For each $S$, let $\Gr(S)$ be the vector space freely generated by $S$-graphs 
and write $\Gr$ for the associated symmetric sequence of vector spaces.
\item If $G\in \Gr(S)$, define $|G|$ to be the cardinality of $S$, which we call the weight of $G$.
 Write $\Gr(n) = \Gr(\{1,...,n\})$.
\end{enumerate}
\end{definition}

We outline the basic properties of the graph cooperad.  For proofs and more  
detailed discussion, see the examples section of \cite{Walt08} where a more 
convenient notation for cooperads is developed.

\begin{definition}
A graph quotient $\phi:G\twoheadrightarrow K$ maps vertices of $G$ to vertices
of $K$ such that edges of $G$ are mapped to either edges of $K$ (with the same
orientation) or vertices of $K$, and
the inverse image of each vertex of $K$ is a non-empty
connected subgraph of $G$. 
\end{definition}

\begin{proposition}\label{P:gr asc}
The symmetric sequence $\Gr$ has a cooperad structure induced by the map 
$$G\longmapsto \!\!\sum_{\phi:G\twoheadrightarrow K} \!\!
          K \ {\textstyle \bigotimes}\ 
          \Bigl(\bigotimes_{k\in \mathrm{Vert}(K)}\phi^\inv(k)\Bigr),$$
where $\phi:G\twoheadrightarrow K$ ranges over all graph quotient maps and 
$\phi^\inv(k)$ is 
the connected subgraph of $G$ mapping to vertex $k$.
\end{proposition}

The cooperad structure above is associative in the sense that the 2-arity 
structure map 
$$\Bigl(\linep{a}{b}\Bigr)^*: G \longmapsto 
\sum_{\phi : G\twoheadrightarrow 
\begin{aligned}\smalllinep{a}{b}\end{aligned}} \phi^\inv(a)\otimes \phi^\inv(b)$$
is (co-)associative.
The symmetric sequence $\Gr$ has another cooperad structure which we call 
anti-commutative for an analogous reason.  

\begin{definition}
Let $E\subset \mathrm{Edge}(K)$.  Define $\mathrm{rev}_E(K)$ to be the 
graph resulting from reversing the orientations of the edges $E$ of $K$.
\end{definition}

\begin{proposition}\label{P:gr acc}
The symmetric sequence $\Gr$ has an anti-commutative cooperad structure induced by 
$$G\longmapsto 
  \!\!\!\sum_{\substack{\phi:G\twoheadrightarrow K \\ E\subset \mathrm{Edge}(K)}}
        \!\!  (-1)^{|E|}\ \mathrm{rev}_E(K) \ {\textstyle \bigotimes}\ 
          \Bigl(\bigotimes_{k\in \mathrm{Vert}(K)}\phi^\inv(k)\Bigr),$$
where $\phi$ and $\phi^\inv(k)$ are as above.
\end{proposition}

\begin{definition}
The anti-commutative graph cooperad, denoted $\grac$, is given by the symmetric sequence
$\Gr$ equipped with the anti-commutative cooperad structure of Proposition~\ref{P:gr acc}.  

The associative graph cooperad, denoted $\gras$, is given by the symmetric
sequence $\Gr$ equipped with the associative cooperad structure of 
Proposition~\ref{P:gr asc}.
\end{definition}



\begin{remark}
Coalgebras over these graph cooperads have not, to our knowledge, been studied before.
We plan to  study  them in future work, but 
there are two main features we would like to highlight now.
First, such graph coalgebras are not binary coalgebras.  For example,
$\Gr(3)$ is twelve-dimensional, while the cooperad structure
map goes to $\Gr(2) \otimes (\Gr(2) \otimes \Gr(1))  \oplus \Gr(2) \otimes (\Gr(2) \otimes \Gr(1))$,
which is eight-dimensional so this structure map cannot be injective.
Secondly, associative graph coalgebras extend associative coalgebras, as 
we establish in \refP{as to gr}.
\end{remark}

In the language of operads, the standard approaches to Lie algebras 
can be summarized by a sequence of operad maps $\Tr \to \Lie \to \Ass$.
Recall that the associative operad 
has $\Ass(n)$ of rank $n!$, naturally spanned by monomials in $n$ variables
with no repetition, and 
$\Tr$ is the the tree operad whose structure maps are defined by grafting
and which governs non-associative binary algebras (see \ref{D:firstdef} below).  
Our Lie coalgebra model follows from fitting the anti-commutative graph cooperad
into the linearly dual sequence of 
cooperads as
$\Tr^{\vee} \leftarrow \Lie^{\vee} \leftarrow \grac \leftarrow \Ass^{\vee}$.  
The following propositions are easily verified by direct calculation.

\begin{proposition}\label{P:as to gr}
The associative cooperad $\Ass^\vee$ maps to the associative graph cooperad $\gras$
by sending the monomial
$x_{1} x_{2} \cdots x_{n}$ to the graph 
$\longgraph{x_1}{x_{2}}{\cdots}{x_{n}}$.
\end{proposition}

\begin{proposition}\label{P:gr to gr}
The associative graph cooperad $\gras$ maps to the anti-commutative
graph cooperad $\grac$ via the map
 $$G \longmapsto \frac{1}{2^{\#(\mathrm{Edge}(G))}} 
 \sum_{E\subset \mathrm{Edge}(G)} (-1)^{|E|} \mathrm{rev}_E(G).$$
\end{proposition}

\medskip 

Next, we develop the configuration pairing between graphs and trees,
which allows us to explicitly 
compute the composition $\grac \to \Tr^\vee$.  We use this to gain a new  
understanding of $\Lie^{\vee}$ ``in the middle.''  
In particular we show in \ref{T:eil dual to lie}  that $\Lie^\vee$ is 
isomorphic as a cooperad to a quotient of $\grac$ which we call $\Eil$.
Furthermore, we show in \ref{T:comp is correct} the standard map 
$\Ass^\vee \to \Lie^\vee$ is equal to 
the composition
of the maps in Propositions~\ref{P:as to gr} and \ref{P:gr to gr} followed by the 
quotient map to $\eil$. 
We first define terms.

\begin{definition}\label{D:firstdef}
Let $S$ be a finite set.
An $S$-tree is an isotopy class of  acyclic 
graphs embedded in the upper half plane 
with all vertices either trivalent or univalent. 
Trivalent vertices are called internal vertices.  
One univalent vertex is distinguished as the root and embedded at the origin. 
The other univalent vertices are called leaves and are equipped with a 
labeling ismorphism $\ell:\mathrm{Leaves}\xrightarrow{\ \cong\ } S$.
We will standardly conflate leaves with their labels.

Let $\Tr(S)$ be the vector space generated by $S$-trees, $\Tr$ be the
associated symmetric sequence of vector spaces, 
and write $\Tr(n)$ for $\Tr(\{1,...,n\})$.
\end{definition}

See  II.1.9 in \cite{MSS02}
for a precise definition of the operad structure maps of $\Tr$ through
grafting.
The pairing between $\Gr(n)$ and $\Tr(n)$ was
developed in \cite{Sinh06.2}, and arises in the study of configuration spaces.  
Let the height of a vertex in a tree be the number of edges between that vertex
and the root.  The nadir of a path in a tree is the vertex of lowest
height which it traverses.

\begin{definition}\label{D:confpair}
Fix a finite set $S$.
Given an $S$-graph $G$ and an $S$-tree $T$, define the map 
$$\beta_{G,T}:\bigl\{\text{edges of } G\bigr\} \longrightarrow 
\bigl\{\text{internal vertices of } T\bigr\}$$ 
by sending an edge from vertex $a$ to $b$ in $G$ to the vertex at the nadir of 
the shortest path in $T$ between the leaves with labels $a$ and  $b$. 
The configuration pairing of
$G$ and $T$ is 
$$\bigl\langle G,\, T\bigr\rangle = 
  \begin{cases} \displaystyle
    \prod_{\substack{e\text{ an edge} \\\text{of }G}} \!\!\!\! 
        \text{sgn}\bigl(\beta_{G,T}(e)\bigr)
                & \text{if $\beta$ is surjective,} \\
         \qquad 0  & \text{otherwise}
\end{cases}$$
where  given an edge $\linep{a}{b}$ of $G$,
$\text{sgn}\left(
  \beta\Bigl(\begin{aligned}\linep{a}{b}\end{aligned}\Bigr)\right) = 1$
   if leaf $a$
is to the left of leaf $b$ under the planar embedding of $T$;  
otherwise it is $-1$.  
\end{definition}

\begin{example}
Following is the map $\beta_{G,T}$ for a single graph $G$ and 
two different trees $T$.
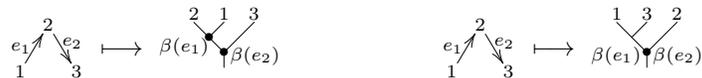
$$\begin{aligned}\begin{xy} 
  (0,-2.5)*+UR{\scriptstyle 1}="a", 
  (3.75,3.75)*+UR{\scriptstyle 2}="b", 
  (7.5,-2.5)*+UR{\scriptstyle 3}="c", 
  "a";"b"**\dir{-}?>*\dir{>} ?(.4)*!RD{\scriptstyle e_1}, 
  "b";"c"**\dir{-}?>*\dir{>} ?(.5)*!LD{\scriptstyle e_2} 
 \end{xy}\end{aligned}\  \longmapsto \ 
 \begin{xy}
   (2,1.5); (4,3.5)**\dir{-}, 
   (0,3.5); (4,-.5)**\dir{-} ?(.5)*!RU{\scriptstyle \beta(e_1)}; 
   (4,-.5); (8,3.5)**\dir{-} ?(.2)*!LU{\scriptstyle \beta(e_2)}, 
   (4,-.5); (4,-2.5)**\dir{-}, 
   (0 ,4.7)*{\scriptstyle 2}, 
   (4,4.7)*{\scriptstyle 1}, 
   (8,4.7)*{\scriptstyle 3}, 
   (2,1.5)*{\scriptstyle \bullet},
   (4,-.5)*{\scriptstyle \bullet},
 \end{xy}  \qquad \qquad \qquad 
 \begin{aligned}\begin{xy} 
  (0,-2.5)*+UR{\scriptstyle 1}="a", 
  (3.75,3.75)*+UR{\scriptstyle 2}="b", 
  (7.5,-2.5)*+UR{\scriptstyle 3}="c", 
  "a";"b"**\dir{-}?>*\dir{>} ?(.4)*!RD{\scriptstyle e_1}, 
  "b";"c"**\dir{-}?>*\dir{>} ?(.5)*!LD{\scriptstyle e_2} 
 \end{xy}\end{aligned}\  \longmapsto \  
 \begin{xy}
   (2,1.5); (4,3.5)**\dir{-}, 
   (0,3.5); (4,-.5)**\dir{-} ?(.8)*!RU{\scriptstyle \beta(e_1)}; 
   (4,-.5); (8,3.5)**\dir{-} ?(.2)*!LU{\scriptstyle \beta(e_2)}, 
   (4,-.5); (4,-2.5)**\dir{-}, 
   (0 ,4.7)*{\scriptstyle 1}, 
   (4,4.7)*{\scriptstyle 3}, 
   (8,4.7)*{\scriptstyle 2}, 
   (4,-.5)*{\scriptstyle \bullet}, 
 \end{xy} $$
In the first example, $\text{sgn}\bigl(\beta(e_1)\bigr) = -1$ and 
$\text{sgn}\bigl(\beta(e_2)\bigr) = 1$. 
In the second example, $\text{sgn}\bigl(\beta(e_1)\bigr) = 1$ and 
$\text{sgn}\bigl(\beta(e_2)\bigr) = -1$. 
The graph and tree of the first example pair to $-1$; in the second example 
they pair to $0$.
\end{example}

From the tree operad, the Lie operad is defined as follows.

\begin{definition}\label{D:Lie}
$\Lie(n)$ is the quotient of $\Tr(n)$ 
by the anti-symmetry and Jacobi relations:

\begin{align*}
\text{(anti-symmetry)} \qquad 
 & \qquad
\begin{xy}
   (0,1.5); (1.5,0)**\dir{-};          
   (3,1.5)**\dir{-},               
   (1.5,-1.5); (1.5,0)**\dir{-},       
   (-.4,2.7)*{\scriptstyle T_1},   
   (3.8,2.7)*{\scriptstyle T_2},   
   (1.5,-2.7)*{\scriptstyle R}
\end{xy}\ =\ - \begin{xy}
   (0,1.5); (1.5,0)**\dir{-};          
   (3,1.5)**\dir{-},               
   (1.5,-1.5); (1.5,0)**\dir{-},       
   (-.4,2.7)*{\scriptstyle T_2},   
   (3.8,2.7)*{\scriptstyle T_1},    
   (1.5,-2.7)*{\scriptstyle R}
\end{xy} \\
\text{(Jacobi)} \qquad 
 & \qquad 
 \begin{xy}   
   (1.5,1.5); (3,3)**\dir{-}, 
   (0,3); (3,0)**\dir{-};
   (6,3)**\dir{-},   
   (3,0); (3,-1.5)**\dir{-}, 
   (-.4,4.2)*{\scriptstyle T_1}, 
   (3.2,4.2)*{\scriptstyle T_2},
   (6.8,4.2)*{\scriptstyle T_3}, 
   (3,-2.7)*{\scriptstyle R}
 \end{xy}  \ + \ \begin{xy}   
   (1.5,1.5); (3,3)**\dir{-}, 
   (0,3); (3,0)**\dir{-};
   (6,3)**\dir{-},   
   (3,0); (3,-1.5)**\dir{-}, 
   (-.4,4.2)*{\scriptstyle T_2}, 
   (3.2,4.2)*{\scriptstyle T_3},
   (6.8,4.2)*{\scriptstyle T_1}, 
   (3,-2.7)*{\scriptstyle R}
 \end{xy} \ + \ \begin{xy}   
   (1.5,1.5); (3,3)**\dir{-}, 
   (0,3); (3,0)**\dir{-};
   (6,3)**\dir{-},   
   (3,0); (3,-1.5)**\dir{-}, 
   (-.4,4.2)*{\scriptstyle T_3}, 
   (3.2,4.2)*{\scriptstyle T_1},
   (6.8,4.2)*{\scriptstyle T_2}, 
   (3,-2.7)*{\scriptstyle R}
 \end{xy} 
 \ =\  0,
\end{align*}
where $R$, $T_1$, $T_2$, and $T_3$ stand for arbitrary (possibly trivial) 
subtrees which are not modified in these operations.
\end{definition}

The configuration pairing respects anti-symmetry and Jacobi relations among 
trees.  There is a similar set of relations which the configuration pairing 
respects among graphs.

\begin{definition}\label{D:eiln}
Let $\Eil(n)$ be the quotient of $\Gr(n)$
by the relations
\begin{align*}
\text{(arrow-reversing)}\qquad & \qquad
\begin{xy}                           
  (0,-2)*+UR{\scriptstyle a}="a",    
  (3,3)*+UR{\scriptstyle b}="b",     
  "a";"b"**\dir{-}?>*\dir{>},         
  (1.5,-5),{\ar@{. }@(l,l)(1.5,6)},
  ?!{"a";"a"+/va(210)/}="a1",
  ?!{"a";"a"+/va(240)/}="a2",
  ?!{"a";"a"+/va(270)/}="a3",
  "a";"a1"**\dir{-},  "a";"a2"**\dir{-},  "a";"a3"**\dir{-},
  (1.5,6),{\ar@{. }@(r,r)(1.5,-5)},
  ?!{"b";"b"+/va(90)/}="b1",
  ?!{"b";"b"+/va(30)/}="b2",
  ?!{"b";"b"+/va(60)/}="b3",
  "b";"b1"**\dir{-},  "b";"b2"**\dir{-},  "b";"b3"**\dir{-},
\end{xy}\ =\ \ -  
\begin{xy}                           
  (0,-2)*+UR{\scriptstyle a}="a",    
  (3,3)*+UR{\scriptstyle b}="b",     
  "a";"b"**\dir{-}?<*\dir{<},         
  (1.5,-5),{\ar@{. }@(l,l)(1.5,6)},
  ?!{"a";"a"+/va(210)/}="a1",
  ?!{"a";"a"+/va(240)/}="a2",
  ?!{"a";"a"+/va(270)/}="a3",
  "a";"a1"**\dir{-},  "a";"a2"**\dir{-},  "a";"a3"**\dir{-},
  (1.5,6),{\ar@{. }@(r,r)(1.5,-5)},
  ?!{"b";"b"+/va(90)/}="b1",
  ?!{"b";"b"+/va(30)/}="b2",
  ?!{"b";"b"+/va(60)/}="b3",
  "b";"b1"**\dir{-},  "b";"b2"**\dir{-},  "b";"b3"**\dir{-},
\end{xy} \\
\text{(Arnold)}\qquad & \qquad
\begin{xy}                           
  (0,-2)*+UR{\scriptstyle a}="a",    
  (3,3)*+UR{\scriptstyle b}="b",   
  (6,-2)*+UR{\scriptstyle c}="c",   
  "a";"b"**\dir{-}?>*\dir{>},         
  "b";"c"**\dir{-}?>*\dir{>},         
  (3,-5),{\ar@{. }@(l,l)(3,6)},
  ?!{"a";"a"+/va(210)/}="a1",
  ?!{"a";"a"+/va(240)/}="a2",
  ?!{"a";"a"+/va(270)/}="a3",
  ?!{"b";"b"+/va(120)/}="b1",
  "a";"a1"**\dir{-},  "a";"a2"**\dir{-},  "a";"a3"**\dir{-},
  "b";"b1"**\dir{-}, "b";(3,6)**\dir{-},
  (3,-5),{\ar@{. }@(r,r)(3,6)},
  ?!{"c";"c"+/va(-90)/}="c1",
  ?!{"c";"c"+/va(-60)/}="c2",
  ?!{"c";"c"+/va(-30)/}="c3",
  ?!{"b";"b"+/va(60)/}="b3",
  "c";"c1"**\dir{-},  "c";"c2"**\dir{-},  "c";"c3"**\dir{-},
  "b";"b3"**\dir{-}, 
\end{xy}\ + \                             
\begin{xy}                           
  (0,-2)*+UR{\scriptstyle a}="a",    
  (3,3)*+UR{\scriptstyle b}="b",   
  (6,-2)*+UR{\scriptstyle c}="c",    
  "b";"c"**\dir{-}?>*\dir{>},         
  "c";"a"**\dir{-}?>*\dir{>},          
  (3,-5),{\ar@{. }@(l,l)(3,6)},
  ?!{"a";"a"+/va(210)/}="a1",
  ?!{"a";"a"+/va(240)/}="a2",
  ?!{"a";"a"+/va(270)/}="a3",
  ?!{"b";"b"+/va(120)/}="b1",
  "a";"a1"**\dir{-},  "a";"a2"**\dir{-},  "a";"a3"**\dir{-},
  "b";"b1"**\dir{-}, "b";(3,6)**\dir{-},
  (3,-5),{\ar@{. }@(r,r)(3,6)},
  ?!{"c";"c"+/va(-90)/}="c1",
  ?!{"c";"c"+/va(-60)/}="c2",
  ?!{"c";"c"+/va(-30)/}="c3",
  ?!{"b";"b"+/va(60)/}="b3",
  "c";"c1"**\dir{-},  "c";"c2"**\dir{-},  "c";"c3"**\dir{-},
  "b";"b3"**\dir{-}, 
\end{xy}\ + \                              
\begin{xy}                           
  (0,-2)*+UR{\scriptstyle a}="a",    
  (3,3)*+UR{\scriptstyle b}="b",   
  (6,-2)*+UR{\scriptstyle c}="c",    
  "a";"b"**\dir{-}?>*\dir{>},         
  "c";"a"**\dir{-}?>*\dir{>},          
  (3,-5),{\ar@{. }@(l,l)(3,6)},
  ?!{"a";"a"+/va(210)/}="a1",
  ?!{"a";"a"+/va(240)/}="a2",
  ?!{"a";"a"+/va(270)/}="a3",
  ?!{"b";"b"+/va(120)/}="b1",
  "a";"a1"**\dir{-},  "a";"a2"**\dir{-},  "a";"a3"**\dir{-},
  "b";"b1"**\dir{-}, "b";(3,6)**\dir{-},
  (3,-5),{\ar@{. }@(r,r)(3,6)},
  ?!{"c";"c"+/va(-90)/}="c1",
  ?!{"c";"c"+/va(-60)/}="c2",
  ?!{"c";"c"+/va(-30)/}="c3",
  ?!{"b";"b"+/va(60)/}="b3",
  "c";"c1"**\dir{-},  "c";"c2"**\dir{-},  "c";"c3"**\dir{-},
  "b";"b3"**\dir{-}, 
\end{xy}\ =\ 0,                             
\end{align*}
where ${a}$, ${b}$, and ${c}$ stand for 
vertices in the graph which could possibly have other connections to
other parts of the graph 
which are not modified in these operations.  We emphasize that $a$, $b$, $c$ 
are {\em vertices}, not subgraphs.  
\end{definition}

Sinha's paper \cite{Sinh06.2} establishes the following theorem, which was first proven
independently by Tourtchine \cite{Tour04} and, in the odd setting, 
Melancon and Reutenauer \cite{MeRe96}.

\begin{theorem}\label{T:operad pairing}
 The configuration pairing 
 $\bigl\langle G,\, T\bigr\rangle$ 
 between $\Gr(n)$ and $\Tr(n)$
 descends to a perfect equivariant pairing between $\eil(n)$ and $\lie(n)$. 

 There is an isomorphism of symmetric sequences $\eil(n) \cong \lie^\vee(n)$.
\end{theorem}

The theorem is proven by first showing that the pairing vanishes
on Jacobi and anti-symmetry combinations of trees as well as on
arrow-reversing and Arnold combinations of graphs.  
These relations allow one to reduce to   
generating sets of ``tall''  trees and 
``long''  graphs -- as in the figure below.  The
pairing is a Kronecker pairing  on these generating sets.

\begin{figure}[ht]
 $$\begin{aligned}\begin{xy}
   (2,1.5); (4,3)**\dir{-},	
   (4,0); (8,3)**\dir{-};
   (0,3); (6,-1.5)**\dir{-};
   (12,3)**\dir{-},
   (6,-1.5);(8,-3)**\dir{.};
   (10,-4.5)**\dir{-};
   (20,3)**\dir{-},
   (10,-4.5);(10,-6)**\dir{-},
   (0,4.2)*{\scriptstyle{1}},
   (4,4.2)*{\scriptstyle{i_2}},
   (8,4.2)*{\scriptstyle{i_3}},
   (12,4.2)*{\scriptstyle{i_4}},
   (20,4.2)*{\scriptstyle{i_n}}
 \end{xy}\end{aligned} \qquad \qquad
 \begin{aligned}\begin{xy}
   (0,-3)*+UR{\scriptstyle 1}="1",
   (4,4.5)*+UR{\scriptstyle j_2}="2",
   (8,-3)*+UR{\scriptstyle j_3}="3",
   (12,4.5)*+UR{\scriptstyle j_4}="4",
   (16,-3)*+UR{\scriptstyle j_{n-1}}="5",
   (20,4.5)*+UR{\scriptstyle j_n}="6",
   "1";"2"**\dir{-}?>*\dir{>},
   "2";"3"**\dir{-}?>*\dir{>},
   "3";"4"**\dir{-}?>*\dir{>},
   "4";"5"**\dir{.},
   "5";"6"**\dir{-}?>*\dir{>},
 \end{xy}\end{aligned}$$
 \caption{Tall trees and long graphs}
\end{figure}
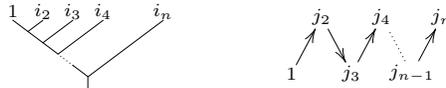


\begin{proposition}\label{P:cooperad coideal}
 The subcomplex of graph expressions generated by 
 arrow-reversing and Arnold expressions of graphs is a coideal \cite[\S 2.1]{GiKa94}
 of $\grac$.
\end{proposition}
 
\begin{corollary}
 The symmetric sequence $\Eil$ inherits an anti-commutative
 cooperad structure from $\grac$.
\end{corollary}

\begin{definition}
By abuse, write $\eil$ for the cooperad induced by quotienting $\grac$ by the
Arnold and arrow-reversing identites.
\end{definition}

\begin{proposition}\label{P:eil dual to lie}
 The cooperad structure of $\grac$
 is compatible with the operad structure of $\Tr$ via the configuration pairing.  
\end{proposition}

\begin{corollary}
 The 
 cooperad structure on $\eil$ is compatible with the operad
 structure of $\lie$ (inherited from that of $\Tr$) via the configuration pairing.
\end{corollary}

\begin{theorem}\label{T:eil dual to lie}
As cooperads, $\eil \cong \lie^\vee$.
Quotienting by Arnold and arrow-reversing identities gives a surjection of
cooperads from $\grac$ to $\lie^\vee$.
\end{theorem}

Since we would rather emphasize free and cofree algebras
than the operads defining them,
we will reserve the computations required for Propositions~\ref{P:cooperad coideal}
and \ref{P:eil dual to lie}
for the the proofs of Propositions~\ref{P:cobracket well defined on E} and 
\ref{P:co/bracket duality} 
which are the analogous statements on the level of coalgebras and algebras.  
A short duality computation (which we leave for the reader)
now completes our operadic picture.

\begin{proposition}\label{T:comp is correct}
The following duality diagram of operads and cooperads commutes.
$$\xymatrix{
 \Lie \ar[rrr] \ar@{<->}[d]^{*} & & & \Ass \ar@{<->}[d]^{*} \\
 \Eil & \grac \ar@{->>}[l] & \gras \ar[l] & \Ass^\vee \ar@{_(->}[l]
}$$
\end{proposition}

Algebra level consequences of this duality are discussed in 
Section~\ref{S:dualPBW} on coenveloping graph coalgebras.

\begin{remark}
Note that our construction of coalgebras is over cooperads 
rather than  
over operads.  It is common in the literature (such as \cite{Smit03})
to largely eschew 
the use of cooperads when discussing coalgebras, instead defining coalgebras over operads 
briefly as follows.
Recall the  endomorphism operad $\mathrm{End}(V)$ of an object $V$ in a closed symmetric monoidal category.  The endomorphism operad of $V$ in 
the opposite category is called its  coendomorphism operad 
$\mathrm{Coend}(V)(n) = \mathrm{Hom}(V,\, V^{\otimes n})$. If $\mathcal{P}$ is an operad
then a $\mathcal{P}$-algebra structure on $V$ is an operad map $\mathcal{P} \to \mathrm{End}(V)$,
and a $\mathcal{P}$-coalgebra structure on $V$ is an operad map $\mathcal{P} \to
\mathrm{Coend}(V)$.

This relates to coalgebras over a cooperad in the following manner.
A map $\mathcal{P} \to \mathrm{Coend}(V)$ consists of equivaraint maps 
$\mathcal{P}(n) \to \mathrm{Hom}(V,\,V^{\otimes n})$.  
If $\mathcal{P}(n)$ is dualizable then these are
the same as equivariant maps $V\to \mathcal{P}(n)^* \otimes V^{\otimes n}$, 
which because $V$ has trivial action are simply maps from $V$ to the
$\Sigma_n$-invariants of the right side.  If $\mathcal{P}(n)$ is dualizable then the 
$\mathcal{P}(n)^*$ form a cooperad, and the structure above defines a coalgebra over this
cooperad.  This  construction is immediately dual to the structure maps 
$\mathcal{P}(n)\otimes V^{\otimes n} \to V$ defining algebras over an operad.
We write $\mathcal{P}^\vee$ for the cooperad 
$\mathcal{P}^\vee(n) = \mathcal{P}(n)^*$.

For more information about a general approach to cooperads and coalgebras over 
cooperads, see \cite{Walt08}.
Developing cooperads on their own terms not only mitigates the use of linear duality,
but gives a more understandable and more computable approach, at least
 to Lie coalgebras and Quillen's rational homotopy theory \cite{Quil69} as we presently
 develop.
\end{remark}

\begin{remark}
While tree operad $\Tr$ governs binary non-associative algebras, the graph cooperads
cannot govern non-associative binary coalgebras.
The configuration pairing between $\Tr$ and 
$\Gr$ is not perfect, nor could there be a different pairing which is perfect.  
For example, 
$\Tr(n)$ has dimension $\frac{n(n-1)}{2} - 1$ as a $\mathbb{Q}[\Sigma_n]$-module
(for $n>1$).  But as a $\mathbb{Q}[\Sigma_3]$-module $\Gr(3)$ is of dimension $3$, and
and as a $\mathbb{Q}[\Sigma_4]$-module $\Gr(4)$ is of dimension $8$.  It is also not
clear what either the linear or Koszul-Moore duals (in the sense of \cite{Moor70})
of graph cooperads are.  
\end{remark}

\section{The pairing between free tree algebras and cofree graph coalgebras}

Constructing our graphical model for Lie coalgebras, we are interested in 
coalgebras over the anti-commutative graph cooperad $\grac$.  Though we may
occasionally write ``anti-commutative graph coalgebra'' for emphasis,
in general we will
write simply ``graph coalgebra'' to mean a coalgebra over the cooperad
$\grac$.  Note that below we explicitly develop only the quadratic structure of 
graph coalgebras since that is all that we require to understand Lie coalgebras. 

\subsection{Basic manipulations of cofree graph coalgebras}

A first step in the theory of operads is the construction of free algebras.
We will use the co-Schur functors associated to $\Gr$ and $\eil$ 
(dual to the Schur functors of \cite{GeJo94}) to construct
explicit models for Lie coalgebras as quotients of anti-commutative graph coalgebras.

\begin{definition}
Let $W$ be a vector space. Define the vector spaces
$\vscofreeG(W)$ and $\vscofreeE(W)$ as follows. 
\begin{equation*}\begin{aligned}
\vscofreeG(W)\ &\cong \ \bigoplus_n\:\left(\Gr(n) \otimes
                            W^{\otimes n} \right)_{\Sigma_n} \\
\vscofreeE(W)\ &\cong \ 
  \bigoplus_n \raisebox{3pt}{$\bigl(\Gr(n) \otimes W^{\otimes n}\bigr)$}
                \!\raisebox{-2pt}{$\diagup$}\!
                \raisebox{-4pt}{$\sim$, $\Sigma_n$}  \ = \ 
  \raisebox{3pt}{$\vscofreeG(W)$}\!\raisebox{-2pt}{$\diagup$}\!\raisebox{-4pt}{$\sim$} 
\end{aligned}
\end{equation*}
where 
$\sim$ is the relation induced by arrow-reversing and Arnold on $\Gr(n)$. 
\end{definition}

There is a difficulty in defining general cofree graph and Lie coalgebras similar 
to that of defining general
cofree associative coalgebras.  Recall that the cotensor coalgebra does not 
give cofree associative coalgebras, since in particular it is always cofinite (that is, 
a finite iteration of the coproduct will reduce any element to primitives). 
Trying to remedy this by replacing colimits by limits usually does
not yield a coalgebra since 
this would require the tensor product to commute with infinite products.
Using results of Smith \cite{Smit03.1}, a cofree graph coalgebra is given in 
general by the largest coalgebra contained in   
${\prod_n}\ \Gr(n)  {\otimes^{\Sigma_n}} W^{\otimes n}$.  Rather than work to get the
correct definition we fall back to the 
time-honored tradition of restricting to $1$-reduced (that is, 
trivial in grading zero and below) coalgebras.  In this category, all coalgebras
are cofinite, the cotensor coalgebra models cofree associative coalgebras, and we have
the following.

\begin{proposition}\label{P:bar E is right}
If $W$ is $1$-reduced, then 
$\vscofreeG(W)$ is the vector space which underlies the cofree graph coalgebra on $W$ and 
$\vscofreeE(V)$ underlies the cofree Lie coalgebra on $V$. 
\end{proposition}


If $W$ is reduced and finitely generated, then so too will
be $\vscofreeG(W)$ and $\vscofreeE(W)$.   We leave the unreduced and infinitely
generated setting for future work.  

We now explicitly develop the graph and Lie coalgebra
structures referred to in the previous proposition.
In the ungraded case, $\vscofreeG(W)$ 
is generated by oriented, connected, acyclic graphs (of possibly infinite size)
whose vertices are labeled by elements of $W$ modulo
multilinearity in the labels.  Cutting a single edge separates graphs in 
$\vscofreeG(W)$, so we may define a coproduct by a summation cutting each edge in turn
and tensoring the resulting graphs in the order determined by the direction of the edge 
which was cut -- this is the coproduct encoded by $\gras$.
In order to descend to the Lie coalgebra cobracket (see Corollary~\ref{C:cobracket is right}) 
we add a twisted term to the above coproduct with signs to make
the result anti-cocommutative -- this is the coproduct encoded by $\grac$.
Explicitly,  $]G[\ =  
\sum_{e \in G} (G^{\hat{e}}_1 \otimes G^{\hat{e}}_2 - 
G^{\hat{e}}_2 \otimes G^{\hat{e}}_1)$, 
where $e$ ranges
over the edges of $G$, and $G^{\hat{e}}_1$ and $G^{\hat{e}}_2$ 
are the connected components
of the graph obtained by removing $e$, which points from 
$G^{\hat{e}}_1$ to $G^{\hat{e}}_2$.

Unfortunately, graded graph coalgebras are more complicated to represent
due to the presence of Koszul signs.
For example, $\graphmp{a}{b}{c}$ could mean either
$\left[\graphmp{2}{1}{3}\bigotimes{b}\otimes{a}\otimes{c}\right]$ or 
$\left[\graphmp{3}{1}{2}\bigotimes{b}\otimes{c}\otimes{a}\right]$,
which differ by a sign of 
$(-1)^{|a||c|}$.  The same difficulty arises when defining graded Lie algebras
via the $\Lie$ operad (or non-associative algebras via the $\Tr$ operad),
but the simple convention there is to choose the equivalence class 
representative whose $\Lie(n)$ component has the ordering of its leaves 
consistent with the planar ordering.  
Because there is no general canonical choice for representativees of  
$\Sigma_n$-equivalence classes in $\Gr(n)$,  we are forced to write  
elements of $\vscofreeG(W)$ explicitly via representatives in 
$\Gr(n) \otimes W^{\otimes n}$.

We define the graded anti-commutative graph cobracket as follows.

\begin{definition}\label{D:cobracket}
The anti-commutative graph cobracket 
$\;]\cdot[\;:\vscofreeG(W) \to \vscofreeG(W) \otimes \vscofreeG(W)$
is given by  
\begin{align*}
\Bigl] G \,{\textstyle \bigotimes}\, w_1 \otimes \cdots \otimes w_n \Bigr[\ = 
 \sum_{e \in G} &(-1)^{\kappa_1} \bigl(G^{\hat{e}}_1 \,{\textstyle \bigotimes}\,
         w_{\sigma^{\hat e}(1)} \otimes \cdots \otimes w_{\sigma^{\hat e}(|G^{\hat e}_1|)}\bigr)
    \bigotimes \bigl(G^{\hat{e}}_2 \,{\textstyle \bigotimes}\,
         w_{\sigma^{\hat e}(|G^{\hat e}_1|+1)}\otimes\cdots\otimes w_{\sigma^{\hat e}(n)}\bigr) \\
 &\ -(-1)^{\kappa_2} \bigl(G^{\hat{e}}_2 \,{\textstyle \bigotimes}\, 
         w_{\sigma^{\hat e}(|G^{\hat e}_1|+1)}\otimes\cdots\otimes w_{\sigma^{\hat e}(n)}\bigr) 
    \bigotimes \bigl(G^{\hat{e}}_1 \,{\textstyle \bigotimes}\,
         w_{\sigma^{\hat e}(1)} \otimes \cdots w_{\sigma^{\hat e}(|G^{\hat e}_1|)}\bigr),
\end{align*} 
where $e$ ranges
over the edges of $G$ and points from the connected subgraph
$G^{\hat{e}}_1$ to the connected subgraph
$G^{\hat{e}}_2$, $\sigma^{\hat e}$ is the unshuffling of vertex labels induced by
separating $G$ into $G^{\hat e}_1$ and $G^{\hat e}_2$, and $(-1)^{\kappa_1},$ 
$(-1)^{\kappa_2}$ are
the Koszul signs due to reordering the $w_i$'s.
\end{definition}

\begin{proposition}
The anti-commutative graph cobracket $\;]\cdot[\;$ on $\vscofreeG(W)$ coincides with the binary 
coproduct arising from the $2$-arity cooperad
structure map of $\grac$.
\end{proposition}

\begin{definition}
Let $\cofreeG(W)$ denote the cofree anti-commutative graph coalgebra on $W$,
whose binary structure is thus given by $\vscofreeG(W)$ with anti-commutative
graph cobracket $\;]\cdot[\;$.  
Similarly, let $\cofreeE(W)$ denote the cofree Lie coalgebra on $W$.
\end{definition}

We will commonly refer to the anti-commutative graph cobracket as merely the
cobracket, since it is the only coproduct operation which we will consider on the 
graph complex $\vscofreeG(W)$.
Our notation will be justified shortly by showing that the anti-commutative graph cobracket 
operation on
$\vscofreeG(W)$ descends to an operation on $\vscofreeE(W)$ 
which coincides with the Lie coalgebra cobracket of $\cofreeE(W)$.  Recall from 
Proposition~\ref{P:bar E is right} that 
$\vscofreeE(W)$ is the vector space underlying $\cofreeE(W)$.

For horizontal brevity we will generally write 
$\overtie{G}{ w_1\otimes \cdots \otimes w_n} := G \bigotimes w_1\otimes \cdots \otimes w_n$ 
for all graphs except for the trivial one: $G = \bullet^1$. 


\begin{example}
The anti-commutative 
graph coalgebra element $\overtie{\graphpm{1}{2}{3}}{a\otimes b\otimes c}$ has cobracket:
\begin{equation*}
\left]\overtie{\graphpm{1}{2}{3}}{a\otimes b\otimes c}\right[\ =
 (-1)^{(|a|+|b|)|c|}\! 
 \left(c\; \otimes
	   \overtie{\linep{1}{2}}{a\otimes b}\right) -
 \left(\overtie{\linep{1}{2}}{a\otimes b}\otimes
	   \; c\right) + 
 \left(a\; \otimes
       \overtie{\linep{2}{1}}{b\otimes c}\right)  -
 (-1)^{|a|(|b|+|c|)}
 \left(\overtie{\linep{2}{1}}{b\otimes c}\otimes
       \; a\right).
\end{equation*}
\end{example}

\begin{proposition}\label{P:cobracket well defined on E}
Let $\mathrm{Arn}(W)$ be the vector subspace of $\vscofreeG(W)$ generated by
arrow-reversing and Arnold expressions of graphs (\ref{D:eiln}).
Then $\mathrm{Arn}(W)$ is a coideal of $\cofreeG(W)$.  That is
$$\bigl]\mathrm{Arn}(W)\bigr[\;\subset \mathrm{Arn}(W)\otimes\cofreeG(W) 
  + \cofreeG(W)\otimes\mathrm{Arn}(W).$$
Thus the cobracket descends to a well-defined operation 
$\,]\cdot[\,:\vscofreeE(W) \to \vscofreeE(W) \otimes \vscofreeE(W)$.
\end{proposition}
\begin{proof}
 Due to the local definition of arrow-reversing and Arnold, it is enough to 
 check the behaviour of the cobracket on an expression reversing the arrow of a 
 graph with only two vertices and on an Arnold expression for a graph with
 only three vertices. 

 The arrow-reversing check (neglecting Koszul signs) is:
$$\left]\overtie{\linep{1}{2}}{a\otimes b} + \overtie{\linem{1}{2}}{a\otimes b} \right[
 \ = \ (a\otimes b - b\otimes a) + (b\otimes a - a\otimes b) \ = \ 0.$$ 

 Modulo arrow-reversing, all graphs with only three vertices are long
 graphs, so it suffices to check the sum of the following (again neglecting signs).
 \begin{align*}
  \left]\overtie{\graphpp{1}{2}{3}}{a\otimes b\otimes c} \right[\ 
    &= 
    \overtie{\linep{1}{2}}{a\otimes b}\otimes \; c\ +\  
    a\; \otimes                                         
         \overtie{\linep{1}{2}}{b\otimes c}\ -\
    c\; \otimes                                         
         \overtie{\linep{1}{2}}{a\otimes b}\ -\ 
    \overtie{\linep{1}{2}}{b\otimes c}\otimes \; a \\   
  \left]\overtie{\graphopp{1}{2}{3}}{a\otimes b\otimes c} \right[\ 
    &= 
    \overtie{\linep{1}{2}}{b\otimes c}\otimes \; a\ +\  
    b\; \otimes                                         
         \overtie{\linep{2}{1}}{a\otimes c}\ -\
    a\; \otimes                                         
         \overtie{\linep{1}{2}}{b\otimes c}\ -\ 
    \overtie{\linep{2}{1}}{a\otimes c}\otimes \; b \\   
   \left]\overtie{\graphpop{1}{2}{3}}{a\otimes b\otimes c}\right[\ 
    &= 
    \overtie{\linep{2}{1}}{a\otimes c}\otimes \; b\ +\  
    c\; \otimes                                         
         \overtie{\linep{1}{2}}{a\otimes b}\ -\
    b\; \otimes                                         
         \overtie{\linep{2}{1}}{a\otimes c}\ -\ 
    \overtie{\linep{1}{2}}{a\otimes b}\otimes \; c \\   
 \end{align*}
\end{proof}

In Proposition~\ref{P:co/bracket duality} below, we show via 
duality that the operation induced
on $\vscofreeE(W)$ by the graph cobracket agrees with the Lie coalgebra cobracket.   In 
Proposition~\ref{P:modcobr} below we prove also the converse of 
Proposition~\ref{P:cobracket well defined on E}:
If $\,]g[\, \in \mathrm{Arn}(W)\otimes\cofreeG(W) 
  + \cofreeG(W)\otimes\mathrm{Arn}(W)$ then $g \in \mathrm{Arn}(W)$.

\begin{remark}\label{R:cobracket}
Though we cannot in general choose canonical representatives of $\Sigma_n$-classes
in $\Gr(n)$, for some classes there is a canonical choice.  
For long $n$-graphs, we chose $\Sigma_n$-representative so that the
ordering of vertices is consistent with the direction of arrows.
In this case we use ``bar'' notation
$$ a_1|a_2|\cdots |a_n  :=
\left[  {\begin{xy}
 (0,-2)*+UR{\scriptstyle 1}="1",
 (3.5,3)*+UR{\scriptstyle 2}="2",
 (7,-2)*+UR{\scriptstyle 3}="3",
 (10.5,3)*+UR{\scriptstyle 4}="4",
 (14,-2)*+UR{\scriptstyle n-1}="5",
 (17.5,3)*+UR{\scriptstyle n}="6",
 "1";"2"**\dir{-}?>*\dir{>},
 "2";"3"**\dir{-}?>*\dir{>},
 "3";"4"**\dir{-}?>*\dir{>},
 "4";"5"**\dir{.},
 "5";"6"**\dir{-}?>*\dir{>},
\end{xy}{\textstyle \bigotimes\,} {a_1\otimes a_2\otimes \cdots \otimes a_n}} \right].$$
Because long $n$-graphs span $\eil(n)$, the bar classes above span
$\vscofreeE(W)$. For example, if $a$, $b$, $c$ and $d$ are all in even degree, then
applying the Arnold and arrow-reversing identites we get
$$\left[ 
\overtie{\begin{xy}
  (0,1)*+UR{\scriptstyle 2}="2",
  "2";"2"+a(90)**{ }?/5mm/*+UR{\scriptstyle 4}="4",
  "2";"2"+a(210)**{ }?/5mm/*+UR{\scriptstyle 1}="1",
  "2";"2"+a(-30)**{ }?/5mm/*+UR{\scriptstyle 3}="3",
  "1";"2"**\dir{-}?>*\dir{>},
  "2";"3"**\dir{-}?>*\dir{>},
  "2";"4"**\dir{-}?>*\dir{>},
\end{xy}}{a\otimes b\otimes c\otimes d}
\right]\ =\ \left[
\overtie{\begin{xy}  
 (0,-2)*+UR{\scriptstyle 4}="1",
 (3.5,3)*+UR{\scriptstyle 2}="2",
 (7,-2)*+UR{\scriptstyle 3}="3",
 (10.5,3)*+UR{\scriptstyle 1}="4",
 "1";"2"**\dir{-}?>*\dir{>},
 "2";"3"**\dir{-}?>*\dir{>},
 "3";"4"**\dir{-}?>*\dir{>},
\end{xy}}{a\otimes b\otimes c\otimes d}
\right]\ -\ \left[  
\overtie{\begin{xy}
 (0,-2)*+UR{\scriptstyle 3}="1",
 (3.5,3)*+UR{\scriptstyle 1}="2",
 (7,-2)*+UR{\scriptstyle 2}="3",
 (10.5,3)*+UR{\scriptstyle 4}="4",
 "1";"2"**\dir{-}?>*\dir{>},
 "2";"3"**\dir{-}?>*\dir{>},
 "3";"4"**\dir{-}?>*\dir{>},
\end{xy}}{a\otimes b\otimes c\otimes d}\right] 
\ =\ 
    d|b|c|a -  
    c|a|b|d.$$
In terms of the bar generators of $\vscofreeE(W)$, the cobracket 
given in 
Proposition~\ref{P:cobracket well defined on E} 
is simply the anti-cocommutative
coproduct (i.e. $\;]\cdot[\ = \Delta - \tau \Delta$ where 
$\tau$ is the twisting map $\tau(x\otimes y) = y\otimes x$). 
This recovers the approach taken by Michaelis \cite{Mich80}
and Schlessinger-Stasheff \cite{ScSt85}.
We elaborate further on this approach in Section~\ref{S:dualPBW}.
\end{remark}

\subsection{Duality of free algebras and cofree coalgebras}

As in the previous section, we start with underlying vector spaces and then move on to
product and coproduct structures.

\begin{lemma}\label{L:Gpair}
Let $G$ be a finite group, and let $V$ and $W$ be modules over a ring in which the order
of $G$ is invertible.  
If $\la-,- \ra$ is an equivariant perfect pairing between $W$ and $V$, 
then the pairing defined between $W_G$ and $V_G$ by 
$\bigl\la [w], [v] \bigr\ra_G = \sum_{g \in G} \la gw, v \ra$ is 
also perfect.
\end{lemma}

\begin{proof}
If $\bigl\la [w], [v] \bigr\ra_G = 0$ for all $[v] \in V_G$ then 
$\bigl\la \sum_{g \in G}  gw, v \bigr\ra = 0$ for all 
$v \in V$.  Because the pairing $\la-,-\ra$ is perfect, this means $ \sum_{g \in G}  gw = 0$ in $W$.
Projecting to $W_G$ implies that $|G|\cdot [w] = 0$, which by our hypotheses means $[w] = 0$.
By equivariance we have $\bigl\la [w], [v] \bigr\ra_G = \sum_{g \in G} \la w, gv \ra$, so we 
may apply the same argument to show that there is no kernel for $\la-,- \ra_G$
in $V_G$ either, yielding the result.
\end{proof}


Let $\freeT(V)$ be the free binary non-associative algebra on $V$, with underlying 
vector space $\vsfreeT(V)$
given by the Schur functor $ \bigoplus_n\:\left(\Tr(n) \otimes
V^{\otimes n} \right)_{\Sigma_n}.$  Define $\freeL(V)$ and $\vsfreeL(V)$ similarly
as the free Lie algebra on $V$ and its underlying vector space.

\begin{definition}
Given $W$ and $V$ vector spaces with a pairing $\la-, -\ra$, 
the configuration pairing between $\vscofreeG(W)$ and 
$\vsfreeT(V)$ is 
$$\bigl\la [G{\textstyle \bigotimes} w_1\otimes\cdots \otimes w_n],\ 
      [T{\textstyle \bigotimes} v_1\otimes\cdots \otimes v_n]\bigr\ra \ = \   
  \sum_{\sigma\in \Sigma_n} \left(
  \bigl\la\sigma G,\, T\bigr\ra \cdot \prod_{i=1}^n
  \la w_{\sigma^\inv(i)},\, v_i\ra\right).$$
\end{definition}

This descends also to a configuration pairing between $\vscofreeE(W)$
and $\vsfreeL(V)$ by Theorem~\ref{T:operad pairing} (proven in 
\cite{Sinh06.2, Tour04, MeRe96}).
Applying \refL{Gpair} we have the following.

\begin{corollary}\label{C:freealgpair}
Over a field of characteristic zero, if $W$ and $V$ pair perfectly then the 
configuration pairing between 
$\vscofreeE(V)$ and $\vsfreeL(W)$ is perfect. 
\end{corollary}
 

\begin{example}
Consider the free Lie algebra on two letters, so that $V$ is spanned
by $a$ and $b$. Then we have the following pairing.  

\begin{align*}
\left\langle 
\left[ \overtie{\graphpp{1}{2}{3}}{a^\ast\otimes a^\ast\otimes b^\ast} \right]
    \ ,\ 
\left[ \overtie{\treel{1}{2}{3}}{a\otimes b\otimes a} \right]
\right\rangle
 \ = \ & 
\left\langle\overtie{\graphpp{1}{2}{3}}{a^\ast\otimes a^\ast\otimes b^\ast}, 
            \overtie{\treel{1}{2}{3}}{a\otimes b\otimes a}\right\rangle \ +\ 
(-1)^{|b||a|}
\left\langle\overtie{\graphpp{1}{3}{2}}{a^\ast\otimes b^\ast\otimes a^\ast}, 
            \overtie{\treel{1}{2}{3}}{a\otimes b\otimes a}\right\rangle \\ &+\ 
(-1)^{|a|^2}
\left\langle\overtie{\graphpp{2}{1}{3}}{a^\ast\otimes a^\ast\otimes b^\ast}, 
            \overtie{\treel{1}{2}{3}}{a\otimes b\otimes a}\right\rangle \ +\ 
\left\langle\overtie{\graphpp{2}{3}{1}}{b^\ast\otimes a^\ast\otimes a^\ast}, 
            \overtie{\treel{1}{2}{3}}{a\otimes b\otimes a}\right\rangle \\ &+\ 
(-1)^{|a|^2+|a||b|}
\left\langle\overtie{\graphpp{3}{1}{2}}{a^\ast\otimes b^\ast\otimes a^\ast}, 
            \overtie{\treel{1}{2}{3}}{a\otimes b\otimes a}\right\rangle \ +\ 
(-1)^{|a|^2}
\left\langle\overtie{\graphpp{3}{2}{1}}{b^\ast\otimes a^\ast\otimes a^\ast}, 
            \overtie{\treel{1}{2}{3}}{a\otimes b\otimes a}\right\rangle  \\
 = \ & 
(-1)^{|b||a|}  
 \left\langle\graphpp{1}{3}{2}, \treel{1}{2}{3}\right\rangle 
\ +\ (-1)^{|a|^2 + |a||b|}
 \left\langle\graphpp{3}{1}{2}, \treel{1}{2}{3}\right\rangle  \\
 =\ &
 0 + (-1)(-1)^{|a|^2+|a||b|}
\end{align*} 
\end{example}

\begin{remark}
 Melan\c{c}on and Reutenauer \cite{MeRe96}
 essentially showed that pairing with with bar elements in $\cofreeG(V^*)$
 defines functionals which can alternately be defined through looking at coefficients
 of Lie polynomials (that is, looking at the coefficients of elements of $\freeL(V)$ in its
 standard embedding in the tensor algebra on $V$).  It would be interesting to understand the functionals
 coming from other elements in $\cofreeG(V^*)$, such as those arising
 from Tourtchine's alternating trees \cite{Tour04}, in a similar manner.
\end{remark}

The configuration pairing further exhibits a duality between non-associative algebra
multiplication and graph 
cobracket operations.  This allows us to compute pairings
inductively.  

\begin{proposition}\label{P:co/bracket duality}
Non-associative algebra multiplication is dual to the anti-commutative graph
cobracket in the configuration pairing.  That is, 
\begin{align*}
 \bigl\langle \gamma,\ (\tau_1\tau_2)\bigr\rangle 
   &= \bigl\langle \,]\gamma[\,,\ \tau_1\otimes \tau_2 \bigr\rangle \\ 
   &= \sum_e \bigl\langle \gamma^{\hat e}_1,\, \tau_1\bigr\rangle\, 
         \bigl\langle \gamma^{\hat e}_2,\, \tau_2\bigr\rangle,
 \end{align*}
where  $]\gamma[\ = \sum_e \bigl(\gamma^{\hat e}_1 \otimes \gamma^{\hat e}_2 \bigr)$. 
\end{proposition}

\begin{proof}
Recall that non-associative algebra multiplication is induced by the 
$\Tr$ operation
$(T_1T_2) = 
\begin{xy}
   (0,1.5); (1.5,0)**\dir{-};          
   (3,1.5)**\dir{-},               
   (1.5,-1.5); (1.5,0)**\dir{-},       
   (-.4,2.7)*{\scriptstyle T_1},   
   (3.8,2.7)*{\scriptstyle T_2},   
\end{xy}$.
We give a bijection between potentially non-zero terms in the 
summands defining $\bigl\la \gamma,\, (\tau_1\tau_2)\bigr\ra$ and
$\bigl\la\,]\gamma[\,,\, \tau_1\otimes \tau_2\bigr\ra$.
In particular, we focus on those terms whose graph/tree pairing component 
may be non-zero. 

Begin by fixing graph and tree representatives.  Let  
$\gamma = [G\otimes \vec{w}] \in 
	  \bigl(\Gr(n)\otimes W^{\otimes n}\bigr)_{\Sigma_n}$ 
(where $\vec{w}\in W^{\otimes n}$) and
$\tau_i = [T_i\otimes \vec{v}_i]\in 
	  \bigl(\Tr(k_i)\otimes V^{\otimes k_i}\bigr)_{\Sigma_{k_i}}$ 
(where $\vec{v}_i \in V^{\otimes k_i}$, $k_1 + k_2 = n$). 
Also let $\gamma^{\hat e}_i = [G^{\hat e}_i \otimes \vec{w}^{\hat e}_i]$ 
(for $i=1,2$)
be the graph coalgebra elements given by 
cutting $\gamma$ at the edge $e$.  
Recall that 
\begin{align*}
 \bigl\langle \gamma,\ (\tau_1\tau_2)\bigr\rangle &= 
  \sum_{\sigma\in \Sigma_n} 
              \bigl\langle \sigma G, (T_1T_2)\bigr\rangle\,
              \bigl\langle \sigma^\inv \vec{w},\, \vec{v}_1\otimes\vec{v}_2 \bigr\rangle\\ 
 \bigl\langle \gamma^{\hat e}_1,\, \tau_1\bigr\rangle\, 
         \bigl\langle \gamma^{\hat e}_2,\, \tau_2\bigr\rangle &=
  \sum_{\sigma_i \in \Sigma_{k_i}}
              \bigl\langle \sigma_1 G^{\hat e}_1,\, T_1\bigr\rangle\,
              \bigl\langle \sigma_2 G^{\hat e}_2,\, T_2\bigr\rangle\,
              \bigl\langle \sigma_1^\inv \vec{w}^{\hat e}_1,\, \vec{v}_1\bigr\rangle\,
              \bigl\langle \sigma_2^\inv \vec{w}^{\hat e}_2,\, \vec{v}_2\bigr\rangle 
 \end{align*}

Suppose that some 
$\bigl\langle \sigma_1 G^{\hat e}_1,\, T_1 \bigr\rangle 
 \bigl\langle \sigma_2 G^{\hat e}_2,\, T_2 \bigr\rangle$ is non-zero. 
Since $G^{\hat e}_1$ and $G^{\hat e}_2$ are the graphs resulting from 
cutting $G$ at the edge $e$,
there is a unique permuation
$\sigma$ which (modulo arrow-reversing at $e$) displays $G$ as
\begin{equation}\label{form2}
 (\sigma G) = \pm
 \begin{xy}
  (0,-3)*+UR{\scriptstyle (\sigma_1 G^{\hat e}_1)}="a",    
  (4,4)*+UR{\scriptstyle (\sigma_2 G^{\hat e}_2) + k_1}="b",     
  "a";"b"**\dir{-}?>*\dir{>} ?(.4)*!RD{\scriptstyle e}       
 \end{xy}
\end{equation}
with sign $\pm$ coming from whether the arrow $e$ was reversed when giving $G$ this form 
(here $(\sigma_2 B_e) + k_1$ denotes adding $k_1$ to each 
vertex label of $(\sigma_2 B_e)$).  
Since the configuration pairing respects the arrow-reversing relation on graphs, 
it follows that
$$\bigl\langle \sigma G,\, (T_1T_2)\bigr\rangle =
\pm \bigl\langle \sigma_1 G^{\hat e}_1,\, T_1\bigr\rangle\,
    \bigl\langle \sigma_2 G^{\hat e}_2,\, T_2\bigr\rangle$$
with the same sign as in Equation~\ref{form2}.

Conversely, if $\bigl\langle \sigma G,\, (T_1T_2)\bigr\rangle$ is non-zero then 
there is a corresponding non-zero 
$\bigl\langle \sigma_1 G^{\hat e}_1,\, T_1\bigr\rangle\,
   \bigl\langle \sigma_2 G^{\hat e}_2,\, T_2\bigr\rangle$.
Given a subset $S \subset \{1,\dots,n\}$ let $G|_S$ denote the full
subgraph of $G$ on the vertices with labels in $S$.
It follows from Definition~\ref{D:confpair}
that $\bigl\langle \sigma G,\, (T_1T_2)\bigr\rangle =0$ unless 
there is exactly one edge in $\sigma G$ between the full subgraphs $(\sigma G)|_{\{1,\dots,k_1\}}$
and $(\sigma G)|_{\{k_1+1,\dots,n\}}$.
Thus these graphs must be connected and (modulo arrow-reversing at $e$) 
the graph $\sigma G$ must be of the form
\begin{equation}\label{form1}
\sigma G = \pm\!\!\!\!\! 
\begin{xy}
  (0,-3)*+UR{\scriptstyle (\sigma G)|_{\{1,\dots,k_1\}}}="a",    
  (4,4)*+UR{\scriptstyle (\sigma G)|_{\{k_1+1,\dots,n\}}}="b",     
  "a";"b"**\dir{-}?>*\dir{>} ?(.4)*!RD{\scriptstyle \sigma e}       
\end{xy}
\end{equation}
with the sign $\pm$ reflecting whether the arrow $\sigma e$ was reversed when writing
$\sigma G$ in this way.
Since the configuration pairing respects the arrow-reversing relation on graphs, it
follows that
$$\bigl\langle \sigma G,\ (T_1T_2)\bigr\rangle =
\pm \bigl\langle (\sigma G)|_{\{1,\dots,k_1\}},\, T_1\bigr\rangle\,
    \bigl\langle (\sigma G)|_{\{k_1+1,\dots,n\}}-k_1,\, T_2\bigr\rangle$$
with the same sign as in Equation~\ref{form1} (where by $G-k$ we mean
to shift all labels of $G$ down by $k$).
We may obtain a non-zero term of the form 
$\langle \sigma_1 G^{\hat e}_1, T_1 \rangle \langle \sigma_2 G^{\hat e}_2, T_2 \rangle$, 
by setting $\sigma_1$ and $\sigma_2$ so that 
$\sigma_1 G^{\hat e}_1 = (\sigma G)_{\{1,\dots,k_1\}}$ and 
$\sigma_2 G^{\hat e}_2 = (\sigma G)_{\{k_1+1,\dots,n\}}-k_1$.

The remainder of the proof is straightforward.  The
signs and pairings between the associated tensors are equal since they
are simply Koszul signs and expected pairings on both sides of the equality.
\end{proof}

The multiplication operation for non-associative algebras 
induces the Lie algebra bracket upon quotienting by
anti-commutativity and Jacobi relations among trees.  In particular, 
Proposition~\ref{P:co/bracket duality} now implies the following.

\begin{corollary}\label{C:cobracket is right}
The 
graph coalgebra cobracket agrees with the Lie coalgebra cobracket through the quotient
map from cofree graph coalgebras to cofree Lie coalgebras.
\end{corollary}

In light of this proposition, an alternate approach to exhibiting the pairing between
$\freeL V$ and $\cofreeG W$  for dual $V$ and $W$ would be 
to define the pairing between $\freeL^n V$ and $\cofreeG^n W$ 
inductively  using the bracket and cobracket.

%

\begin{remark}
 Corollary~\ref{C:freealgpair} and Proposition~\ref{P:co/bracket duality} give a method for 
 constructing functionals on Lie algebras which are not free.
 Any finitely generated graded
 Lie algebra is the homology of some free finitely generated differential graded Lie algebra.
 That is, $L \cong H_* (\freeL V,\, d)$.  The complex $(\freeL V,\, d)$ is dual to 
 $(\cofreeE V^*\!,\, d^*)$,
 whose homology pairs with that of $H_* (\freeL V,\, d)$, namely $L$, 
 through the configuration pairing.  Using bar
 basis elements from $\cofreeE V^*$ one can recover the embedding
 of $L$ in its universal enveloping algebra, but the approach through $\cofreeE V^*$
 offers more flexibility.
\end{remark}

\subsection{Coenveloping graph coalgebras}\label{S:dualPBW}

There are four basic approaches to the free Lie algebra $\freeL(V)$ on a vector space
$V$.
\begin{enumerate}
\item 
 $\freeL$ is  the left adjoint of the forgetful functor from Lie algebras to
 vector spaces.
\item 
 $\freeL(V) \cong \bigoplus_n \lie(n) \otimes_{\Sigma_n} V^{\otimes n}$,
 where $\sigma \in \Sigma_n$ acts on $\lie(n) \otimes V^{\otimes n}$ as
 $\sigma \otimes \sigma^{-1}$, and the $\Sigma_n$ action on $V^{\otimes n}$
 is governed by the Koszul sign convention. 
 \item 
 $\freeL(V)$ is a quotient of the free non-associative algebra on $V$,
 $\bigoplus_n (\Tr(n) \otimes V^{\otimes n})_{\Sigma_n}$,   
 by the anti-symmetry and Jacobi relations on $\Tr(n)$.
\item \label{l3}
 $\freeL(V)$ is the smallest subspace of the tensor algebra on $V$ 
 which contains $V$ and is closed under commutators. 
\end{enumerate}

So far our development of Lie coalgebras has paralleled the second and third
approaches, while the adjointness properties are immediate.  To complete
our picture, and connect with previous work, we
now focus on developing the last approach.
We give a representation of $\cofreeE(W)$
which is dual to  the Poincar\'e-Birkhoff-Witt
embedding of $\freeL(V)$ in the tensor algebra $TV$.
We will exhibit $\cofreeE(W)$ as a quotient of
the cotensor coalgebra. 
This representation is the starting point for the seminal work of Michaelis \cite{Mich80}
on Lie coalgebras, so we in particular identify how our graph model for cofree Lie coalgebras
encompasses that approach.   

\begin{definition}
 Define the graded vector space
 $\modcobr{\vscofreeG(W)}\,$ inductively, setting $\modcobr{\vscofreeG^1(W)}\, = W$ and
 letting $\modcobr{\vscofreeG^n(W)}$ be the quotient of $\vscofreeG^n(W)$ by
 the kernel of the map 
 $$\,]\cdot[\,:\vscofreeG^n(W)\longrightarrow\left(\modcobr{\vscofreeG^{<n}(W)}\right) 
 {\textstyle \bigotimes} \left(\modcobr{\vscofreeG^{<n}(W)}\right).$$
\end{definition}

\begin{proposition}\label{P:modcobr}
 $\vscofreeE(W) \cong \modcobr{\vscofreeG(W)}\,$.  
\end{proposition}

We encourage the interested reader to work through a direct proof of this proposition
by explicitly showing the converse of Proposition~\ref{P:cobracket well defined on E}.
Instead we use duality and compatibility of graph cobrackets with Lie brackets.

\begin{proof}
 By Proposition~\ref{P:cobracket well defined on E}, 
 $\mathrm{ker}(]\!\!\cdot\!\![)\,\supset \mathrm{Arn}(W)$.  It
 remains to show only that $\mathrm{ker}(]\!\!\cdot\!\![)\, \subset \mathrm{Arn}(W)$.
 By Corollary~\ref{C:freealgpair} it is enough to show that the kernel of the pairing between
 $\cofreeG^n(W)$ and $\freeL^n(V)$ contains $\mathrm{ker}(]\!\cdot\![)$. 
 This follows by induction using Proposition~\ref{P:co/bracket duality}. 
\end{proof}

Proposition~\ref{P:modcobr} implies that $\cofreeE(W)$ is the quotient
of $\cofreeG(W)$ by the largest coideal in the kernel of $\cofreeG(W)\to W$.
This extends the definition of Lie coalgebras given by \cite{ScSt85} 
as the quotient of the cotensor coalgebra $TW$ by the largest coideal 
in the kernel of $TW\to W$.  In particular the construction of \cite{ScSt85}
follows as an immediate corollary
using the injection of operads $\Ass^\vee \to \gras$.
We record this in a more computationally useful form as follows.

\begin{corollary}\label{P:T modcobr}
$\vscofreeE(W)$ is isomorphic to the quotient of the cotensor coalgebra $TW$ by the
non-primitive kernel of 
the anti-cocommutative coproduct.
\end{corollary}

\begin{proof}
 There is a ``graphification'' map $g$ which
 injects $TW$ into $\cofreeG W$:
 $$g:w_1|w_2|\cdots|w_n \longmapsto
 \left[\overtie{\begin{xy}
 (0,-2)*+UR{\scriptstyle 1}="1",
 (3.5,3)*+UR{\scriptstyle 2}="2",
 (7,-2)*+UR{\scriptstyle 3}="3",
 (10.5,3)*+UR{\scriptstyle 4}="4",
 (14,-2)*+UR{\scriptstyle n-1}="5",
 (17.5,3)*+UR{\scriptstyle n}="6",
 "1";"2"**\dir{-}?>*\dir{>},
 "2";"3"**\dir{-}?>*\dir{>},
 "3";"4"**\dir{-}?>*\dir{>},
 "4";"5"**\dir{.},
 "5";"6"**\dir{-}?>*\dir{>},
 \end{xy}}{w_1\otimes w_2\otimes \cdots \otimes w_n}\right].$$ 
 By abuse, call the anti-cocommutative coproduct on $TW$ the cobracket,
 and denote it $\;]\cdot [\; = \Delta - \tau\, \Delta$.
 As mentioned in Remark~\ref{R:cobracket}, $g$ 
 sends cobrackets of cotensors in $TW$ to cobrackets of long graphs 
 in $\cofreeE W$.
Now apply Proposition~\ref{P:modcobr}.
\end{proof}


Proposition~\ref{P:modcobr} suggests a simple algorithm for checking whether a Lie
coalgebra element is trivial.  Inductively define the iterated cobracket on graph 
coalgebras
$]\cdot[^n\,:G \to G^{\otimes n}$ by
$$
]g[^n\ = \sum_e \ ]g^{\hat e}_1[^{n-1}\,\otimes\, g^{\hat e}_2,
$$
where $]g[\ = \sum_e g^{\hat e}_1 \otimes g^{\hat e}_2$.
By Proposition~\ref{P:modcobr}, a necessary condition for a graph expression 
$g \in \cofreeE^n(W)$ to 
be trivial is for $]g[^{n-1}\; = 0$.  In fact, this condition is also 
sufficient.

\begin{proposition}\label{P:iterated cobr}
 An element $g\in \cofreeE^n(W)$ is trivial if and only if $]g[^{n-1}\; = 0$.
\end{proposition}
\begin{proof}
 Applying Proposition~\ref{P:co/bracket duality},
 $$\left\la g,\ \Bigl[ \bigl[ [v_1,v_2],v_3\bigr],\cdots v_n\Bigr]\right\ra \ = \ 
  \Bigl\la ]g[^{n-1}\,,\ v_1\otimes v_2 \otimes \cdots \otimes v_n\Bigr\ra.$$
 Since bracket expressions of the form $\bigl[[[v_1,v_2],v_3]\cdots v_n\bigr]$
 span $\freeL W^*$ and the configuration pairing is perfect between 
 $\cofreeE W$ and $\freeL W^*$, $g=0$ if and only if $]g[^{n-1}\;= 0$.
\end{proof}

Our recovery of the approaches to cofree Lie coalgebras of
Michaelis  \cite{Mich80} and Schlessinger-Stasheff \cite{ScSt85}
allows us to highlight some advantages of the graph model.  
Working from $\cofreeG(W)$ the list of relations satisfied by Lie coalgebra elements
is relatively simple to describe -- arrow-reversing and Arnold relations for graphs
as well as symmetric group action.
Once we have restricted to the bar generators, however, the relations become  harder to
describe.  
For example, below are a two relations satisfied by bar generators of
$\cofreeE^n W$ (neglecting Koszul signs). 
 \begin{align} 
     \bigl(w_1|w_2|\cdots|w_n\bigr) -                           
                  (-1)^{n-1}\bigl(w_n|\cdots|w_2|w_1\bigr) = 0 \label{reverse all} \\
     \sum_{\substack{\sigma\text{ a cyclic}\\ \text{permuation of }(1,\ldots,n)}} \kern -20pt
                  \bigl(w_{\sigma(1)}|\cdots|w_{\sigma(n)}\bigr) = 0\label{cyclic rotation}
 \end{align}

Relation (\ref{reverse all}) above comes from applying the arrow-reversing identity
at every arrow of a long graph.
Relation (\ref{cyclic rotation}) is easily verified using 
Proposition~\ref{P:iterated cobr}.  To complete the comparison to 
\cite{ScSt85} we use our graph model to show that quotienting cotensor coalgebras
by shuffle relations gives Lie coalgebras. 

\begin{proposition} \label{Harrison}
 The Harrison shuffles
 give a spanning set of relations among bar generators of
 $\cofreeE^n W$; i.e. 
 $$
     \sum_{\substack{\sigma\text{ a shuffle of} \\
                      (1,2,\dots k)\text{ into }(k+1,\dots,n)}} \kern -27pt
         \bigl(w_{\sigma(1)}|\cdots|w_{\sigma(n)}\bigr) 
           = 0.
 $$
\end{proposition}
\begin{proof}
 Write $\mathrm{Sh}(W)$ for the vector subspace of $\cofreeG W$ generated by the
 Harrison shuffles of bar expressions.  
 It is straightforward to show that $\mathrm{Sh}(W)$ is a coideal:
$$\bigl]\mathrm{Sh}(W)\bigr[\ \subset \; \mathrm{Sh}(W)\otimes\cofreeG W +
                                \cofreeG W\otimes\mathrm{Sh}(W).$$ 
 On bar expressions of either 2 or 3 elements, the Harrison shuffles are merely the 
 arrow-reversing and Arnold 
 relations. Thus by Proposition~\ref{P:iterated cobr}, $\mathrm{Sh}(W) \subset \mathrm{Arn}(W)$.

 That $\mathrm{Sh}(W)$ gives all relations among bar generators is now an immediate application
 Proposition~\ref{P:modcobr} and comments at the end of the first section of \cite{ScSt85}.
\end{proof}

Note that it is not at all clear that relations (\ref{reverse all}) and (\ref{cyclic rotation})
above are inside the coideal of Harrison shuffles.
For computational purposes, it is convenient to have a more minimal set of relations among 
bar generators of $\cofreeE W$.  
Directly applying the configuration pairing, we find the following set of relations.

\begin{proposition} \label{lastReln}
 The below shuffles give a spanning set of relations among bar generators of 
 $\cofreeE^n W$. 
 $$
     \bigl(w_1|w_2|w_3|\cdots|w_n\bigr) + \;(-1)^k \kern -40pt
         \sum_{\substack{\sigma\text{ a shuffle of} \\
                      (k-1,k-2,\dots 1)\text{ into }(k+1,\dots,n)}} \kern -37pt
         \bigl(w_k|w_{\sigma(1)}|\cdots|w_{\sigma(n-1)}\bigr)
           = 0
 $$
\end{proposition}
\begin{proof}
 Let $V$ have basis $v_1,\dots,v_n$ dual to $w_1,\dots,w_n$.  Recall that 
 $\freeL^nV$ is generated
 by Lie bracket expressions $\bigl[[[v_1, v_{i_1}], v_{i_2}],\dots v_{i_{n-1}}\bigr]$.  
 From the definition of the configuration pairing, it follows that the long graph
 $g = w_{j_1}|\cdots|w_{j_m}|w_1|w_{j_{m+1}}|\cdots|w_{j_{n-1}}$ pairs nontrivially 
 with a generating Lie
 bracket expression if and only if $({i_1}, {i_2}, \dots, {i_{n-1}})$ is a shuffle
 of $({j_m}, \dots, {j_1})$ into $({j_{m+1}}, \dots, {j_{n-1}})$ and in this case
 pairs to $(-1)^m$.  We may thus express $g$ in terms of the dual generating
 long graphs $w_1|w_{i_1}|\cdots|w_{i_{n-1}}$ as 
\begin{equation}\label{eq1}
w_{j_1}|\cdots|w_{j_m}|w_1|w_{j_{m+1}}|\cdots|w_{j_{n-1}} = 
 (-1)^{m} \kern -40pt 
    \sum_{\substack{\sigma\text{ a shuffle of} \\
          (j_m,\dots, j_1)\text{ into }(j_{m+1},\dots,j_{n-1})}} \kern -37pt
        \bigl(w_1|w_{\sigma(1)}|\cdots|w_{\sigma(n-1)}\bigr).
\end{equation}
 This is a complete set of relations since it expresses every long graph in terms
 of generating elements.  Relettering so that 
 $w_{j_1}|\cdots|w_{j_m}|w_1|w_{j_{m+1}}|\cdots|w_{j_{n-1}}$ becomes 
 $w_1|\cdots|w_n$ we have the desired relations.
\end{proof}

Equation~\ref{eq1} may be of independent interest, since it gives rise to a canonical 
vector space basis for
cofree Lie coalgebras.  To our knowledge, bases of free Lie algebras involve making choices.

\section{Bar constructions to and from the category of graph coalgebras}

Throughout this section we use $\cofreeC V$
to mean the cofree graded-cocommutative coalgebra on a vector space $V$. 
If $V$ is reduced then $\cofreeC V$ is given by the symmetric invariants of the 
cotensor coalgebra $T^cV$ on $V$ (where the symmetric group acts with Koszul signs). 
Working rationally (with $V$ finitely generated), the norm map gives a vector
space isomorphism with $\freeA V$, 
the free graded-commutative algebra generated by $V$, which is 
given by the symmetric coinvariants of the tensor algebra $TV$ on $V$. 
Elsewhere in the literature this is sometimes called $\Lambda V$ or $S V$.
Our notation is inspired by the standard notation of $\freeL V$ for 
the free Lie algebra on $V$ as well as our mirroring notation $\cofreeG W$ for
the cofree graph coalgebra on $W$.  

Note that $\cofreeC^0 V = \langle 1 \rangle = \freeA^0 V,$ while 
$\freeL^0 V = 0 =\cofreeE^0 V$.  In various instances  we will  take 
augmentation ideals of algebras (denoted $\bar A$) or coaugmentation coideals of 
coalgebras (denoted $\bar C$).

\subsection{The Quillen functors $\eL$ and $\C$} 
  
Recall the standard definition of the Quillen adjoint pair of functors 
$\eL:\mathrm{\dgc} \rightleftarrows \mathrm{\dgl}:\C$. 
The functor 
$\eL$ can be viewed as the cobar construction 
followed by taking Hopf algebra primitives;
$\C$ can be viewed as the bar construction on the universal enveloping algebra 
of a Lie algebra.
Topologically these are identifying the rational homotopy of a space inside the 
cohomology of its loopspace via the Milnor-Moore theorem.
In explicit algebra, given a differential graded-cocommutative 
coalgebra $(C,\ \Delta_C,\ d_C)$,  the functor $\eL$ 
produces the free graded Lie algebra on $s^\inv\bar C$ with a differential
consisting of the free extension of the differential $d_C$ plus a  
``twisting differential'' freely induced by $\Delta_C$.  
Explicitly, we have the following.

\begin{definition}\label{D:eL}
Let $\eL : \dgc \to \dgl$ be the total complex of the bicomplex 
$$\eL(C,\ \Delta_C,\ d_C) = 
\bigl(\freeL(s^\inv\bar C),\  
    d_{\freeL s^\inv C},\  d_\Delta\bigr),$$ 
where $d_{\freeL s^\inv C}$ is the differential inherited from the
differential $d_C$ on $C$; and $d_\Delta$ is the free extension of the 
map given on the generators 
of $\freeL(s^\inv \bar C)$ by 
$$d_\Delta(s^\inv c) = 
{\textstyle \frac{1}{2}}\sum_i (-1)^{|a_i|}[s^\inv a_i, s^\inv b_i]\qquad 
   \mathrm{where}\ \bar\Delta_C c = \sum_i a_i \otimes b_i.$$
\end{definition}

The functor $\C$ is defined dually -- given $(L,\ [\cdot,\cdot]_L,\ d_L)$ in $\dgl$, 
the functor $\C$ takes this to the 
cofree graded-cocommutative coalgebra primitively cogenerated by $sL$ with a differential
consisting of the cofree extension of the differential $d_L$ plus a 
``twisting differential'' cofreely induced by the bracket $[\cdot,\cdot]_L$.  

\begin{definition}
Let $\C : \dgl \to \dgc$ be the total complex of the bicomplex
$$\C(L,\ [\cdot,\cdot]_L,\ d_L) = 
 \bigl(\cofreeC (s L),\ 
    d_{\cofreeC s L},\  d_{[\cdot,\cdot]}\bigr),$$
where $d_{\cofreeC s L}$ is the differential inherited from $d_L$ on $L$; 
and $$d_{[\cdot,\cdot]}(sv_1\cdot sv_2 \cdots sv_n) =
\sum_{i<j} (-1)^{n_{ij}+|v_i|} s[v_i, v_j] \cdot
 sv_1 \cdots \widehat{sv_i} \cdots \widehat{sv_j} \cdots sv_n,$$
where $(-1)^{n_{ij}}$ is the Koszul sign change incurred by moving $sv_i$, and $sv_{j}$ 
to the beginning of this expression.
\end{definition}

An alternate way to view $d_{[\cdot,\cdot]}$, in parallel to \refD{eL}, is
the following.

\begin{proposition}
The differential $d_{[\cdot,\cdot]}$ is the cofree extension of the graded
vector space map $[\cofreeC sL]_\g \to s[L]_\g$ given on 
on $\cofreeC^{\neq 2} sL$ by the zero map and on $\cofreeC^2 sL$ by  
$(sv_1\cdot sv_2)  \mapsto
  (-1)^{|v_1|}s[v_1,v_2].$
\end{proposition}

Adjointness of $\eL$ and $\C$ follows from that of the bar
and bar construction as well as that of the universal enveloping algebra
and Lie primitives functors.  

\begin{remark}\label{R:hat L}
We will shortly construct functors $\E:\dga \leftrightarrows \dge: \A$ (dual to 
$\eL$ and $\C$) as quotients of functors $\G:\dga \leftrightarrows \dgg: \hat \A$
to and from graph coalgebras.  

We may attempt to define a functor $\hat \eL:\dgc \to \dgt$ (where $\dgt$ denotes
dg-non-associative binary algebras and $\mathbb{T}$ denotes the free such algebra) by 
$$
\hat \eL(C,\;\Delta_C,\;d_C) = \bigl(\mathbb{T}(s^\inv \bar C),\;d_{\mathbb{T}s^\inv\bar C},\;
                                  d_\Delta\bigr). 
$$
Unfortunately, $d_\Delta$ is not a differential on the non-associative algebra
$\mathbb{T}(s^\inv \bar C)$ so $\hat \eL$ isn't a differential complex.   
This is a striking difference between the non-associative
algebra approach to Lie algebras and the graph coalgebra approach to Lie coalgebras
we present in the next section.  Indeed, we could presumably
replace non-associative algebras by graph algebras in the above construction and get 
a functor $\hat \eL$ which mapped to differential graded
complexes and generalized both Adams' bar construction and Quillen's $\eL$ 
functor appropriately.  We leave that for furture work.
\end{remark}

\subsection{The functor $\G$}\label{S:E}

To define $\G$, 
we start with a differential graded (commutative, augmented, unital) algebra 
$(A,\ \mu_A,\ d_A)$, with augmentation ideal $\bar A$.  The functor
$\G$ produces the cofree graded anti-commutative graph coalgebra 
on $s^\inv\bar A$ with differential consisting of the 
cofree extension of
the differential $d_A$ along with another part coming from the 
multiplication $\mu_A$, defined by contracting edges.  In order to make this precise, 
we must carefully define the sign associated to contracting an edge.

\begin{definition}
Let $[g]$ be a homogeneous element of $\cofreeG(s^\inv {\bar A})$, namely
an ordered directed graph with 
$n$ vertices along with a tensor of $n$
elements of $s^\inv\bar A$ modulo the usual $\Sigma_n$-action.  
For every edge $e$ of $g$ we may construct a new ordered labeled graph $\mu_e(g)$
as follows.

Pick a representative of $[g]$ modulo $\Sigma_n$
in which edge $e$ goes from 
vertex number $1$ to vertex number $2$, with the first two entries
of the associated tensor being $a$ and $b$.
Contract the edge from $1$ to $2$ in this representative to a vertex 
which is then given the number
$1$ and first entry in the tensor of  $(-1)^{|a|}s^\inv (a b)$.
In this operation, the ordering of all other vertices in the graph is 
shifted down by one to make up for the now missing 2 (associated elements
in the tensor remain the same).
$$\mu_e:
\overtie{
\begin{xy}                           
  (0,-2)*+UR{\scriptstyle 1}="a",    
  (3,3)*+UR{\scriptstyle 2}="b",     
  "a";"b"**\dir{-}?>*\dir{>} ?(.4)*!RD{\scriptstyle e\,},         
  (1.5,-5),{\ar@{. }@(l,l)(1.5,6)},
  ?!{"a";"a"+/va(210)/}="a1",
  ?!{"a";"a"+/va(240)/}="a2",
  ?!{"a";"a"+/va(270)/}="a3",
  "a";"a1"**\dir{-},  "a";"a2"**\dir{-},  "a";"a3"**\dir{-},
  (1.5,6),{\ar@{. }@(r,r)(1.5,-5)},
  ?!{"b";"b"+/va(90)/}="b1",
  ?!{"b";"b"+/va(30)/}="b2",
  ?!{"b";"b"+/va(60)/}="b3",
  "b";"b1"**\dir{-},  "b";"b2"**\dir{-},  "b";"b3"**\dir{-},
\end{xy}
}{s^\inv a\otimes s^\inv b\otimes\cdots}\ \longmapsto\ 
(-1)^{|a|}\!\!\!\!\!\!\!\overtie{
\begin{xy}
  (0,1)*+UR{\scriptstyle 1}="a",
  (0,-3),{\ar@{. }@(l,l)(0,4)},
  ?!{"a";"a"+/va(210)/}="a1",
  ?!{"a";"a"+/va(240)/}="a2",
  ?!{"a";"a"+/va(270)/}="a3",
  "a";"a1"**\dir{-},  "a";"a2"**\dir{-},  "a";"a3"**\dir{-},
  (0,4),{\ar@{. }@(r,r)(0,-3)},
  ?!{"a";"a"+/va(90)/}="b1",
  ?!{"a";"a"+/va(30)/}="b2",
  ?!{"a";"a"+/va(60)/}="b3",
  "a";"b1"**\dir{-},  "a";"b2"**\dir{-},  "a";"b3"**\dir{-},
\end{xy}}{s^\inv (ab)\otimes \cdots}
$$
\end{definition}

\begin{definition}\label{D:E}
Let $\G : \dga \to \dgg$ be the total complex of the bicomplex
$$\G(A,\ \mu_A,\ d_A) =
 \bigl(\cofreeG(s^\inv\bar A),\ 
  d_{\cofreeG s^\inv\!\bar A},\ d_\mu\bigr),$$
where $d_{\cofreeG s^\inv\!\bar A}$ takes 
$d_{s^\inv\!\bar A} = -s^\inv d_A$ term-wise in the tensor associated
to a graph coalgebra element and
$d_\mu ([g]) = 
  \sum_{e} [\mu_e g].$
\end{definition}

This is a bicomplex by the same calculation which shows  
that Adams' classical bar construction is a bicomplex.  Indeed, 
$\G$ extends Adams' bar construction to the
category of graph coalgebras.

\begin{proposition}\label{P:dmu}
The map $d_\mu$ is compatible with the cobracket on $\cofreeG(s^\inv\bar A)$.
That is, 
$d_\mu(\,]g[\,) =\ ]d_\mu g[$.

Moreover, $d_\mu$ is the 
cofree extension of the graded vector space map 
$\bigl[\cofreeG s^\inv \bar A\bigr]_\g \to s^\inv [\bar A]_\g$ given 
on $\cofreeG^{\neq 2} s^\inv \bar A$ by the zero map and on  
on $\cofreeG^2 s^\inv \bar A$
by $\linep{1}{2}\!{\textstyle \bigotimes} s^\inv a\otimes s^\inv b\ \longmapsto\ 
 (-1)^{|a|} s^\inv (ab).$
\end{proposition}


We now construct our Lie coalgebraic bar construction $\E$
as a quotient of our graphical bar constuction $\G$.

\begin{proposition}\label{P:d_mu well-def}
 The differential $d_\mu$ preserves the vector subspace generated by arrow-reversing and 
 Arnold expressions.
 Thus the arrow-reversing and Arnold coideal is a subcomplex of $\G$.
\end{proposition}

\begin{proof}
This proposition follows immediately from the compatibility of 
$d_\mu$ and the cobracket and Proposition~\ref{P:cobracket well defined on E}
once we show that $d_\mu$ vanishes on arrow-reversing expressions.
Using the bar representation for graphs, this is shown by:
$$d_\mu\Bigl(s^\inv a| s^\inv b + (-1)^{(|a|+1)(|b|+1)} s^\inv b| s^\inv a\Bigr) \ = \
(-1)^{|a|}s^\inv (ab) + (-1)^{|a| |b|+|a|+1} s^\inv(ba) \ = \ 0.$$
\end{proof}

\begin{definition}
 Let $\E(A)$ be $\G(A)$ modulo the arrow-reversing and Arnold subcomplex.  
\end{definition}

\begin{remark}
In terms of the bar generators, the differentials in the definition of $\E$ 
coincide with the differentials used to define the usual algebraic
(associative) bar construction, but are now defined on the quotient of
the bar construction by the relations induced by arrow-reversing and Arnold.
By Proposition~\ref{Harrison}, $\E(A)$ is isomorphic to the Harrison complex  
of the commutative algebra $A$ equipped with the Lie coalgebra structure from
\cite{ScSt85}.
\end{remark}

\subsection{The functor $\hat \A$}\label{S:A}

The functor $\hat \A$ is given by Adams' cobar construction applied to a 
graph coalgebra.
Explicitly it takes the differential graded graph coalgebra
$(G,\ ]\cdot[_G,\ d_G)$
to the free graded-commutative algebra generated by $sG$ with a differential consisting
of the free extension of $d_G$ along with another part
coming from the graph cobracket.  

\begin{definition}\label{D:A}
Let $\hat \A : \dgg \to \dga$ be the total complex of the bicomplex
$$\hat \A(G,\ ]\cdot[_G,\ d_G) =
\bigl(\freeA s G,\ d_{\freeA s G},\ d_{]\cdot[_G}\bigr),$$
where $d_{]\cdot[}$ is the free extension of the map given on the 
generators of $\freeA sG$ by 
$$d_{]\cdot[}(sg) = 
 {\textstyle \frac{1}{2}}
  \sum_e (-1)^{|g^{\hat e}_1|}\,sg^{\hat e}_1\cdot sg^{\hat e}_2,\qquad  \text{for\ \ } 
   ]g[\; = \sum_e g^{\hat e}_1\otimes g^{\hat e}_2.$$
\end{definition}

Unlike in Remark~\ref{R:hat L}, this defines a differential graded complex. 

\begin{theorem}\label{T:hat A is good}
 $\hat \A(G)$ is a bicomplex.
\end{theorem}
\begin{proof}
 We already know $d_{\freeA sG}^2 = 0$.  
 Also, that $d_{]\cdot[} d_{\freeA sG} = d_{\freeA sG} d_{]\cdot[} = 0$
 follows from anti-cocommutativity of the cobracket.

 To show $d_{]\cdot[}^2 = 0$,
 it is enough to show that $d_{]\cdot[}^2 = 0$ on $sG \subset \freeA sG$. 
 Furthermore it is enough to show $d_{]\cdot[}^2$ vanishes on graphs with only
 three vertices, since the general case is then solved by replacing vertices
 by graphs.  
 $$\begin{aligned}
   d_{]\cdot[}^2 \left(s\overtie{\graphpp{1}{2}{3}}{a\otimes b\otimes c}\right) &=
	d_{]\cdot[}\left(
	   (-1)^{|a|+|b|}s\Bigl(\overtie{\linep{1}{2}}{a\otimes b}\Bigr) \cdot s c  +
	   (-1)^{|a|}s a \cdot s\Bigl(\overtie{\linep{1}{2}}{b\otimes c}\Bigr)\right) \\
   &= (-1)^{|a|+|b|+|a|} sa \cdot sb \cdot sc + 
	  (-1)^{|a|+(|a|+1)+|b|} sa \cdot sb \cdot sc \ = \ 0 
 \end{aligned}$$
 The computations for $\Bigl(\overtie{\graphmp{1}{2}{3}}{a\otimes b\otimes c}\Bigr)$ and 
 $\Bigl(\overtie{\graphpm{1}{2}{3}}{a\otimes b\otimes c}\Bigr)$ are similar (though the
 signs involved are slightly more unpleasant).
\end{proof}

The following proposition is an immediate consequence of
Proposition~\ref{P:cobracket well defined on E}. 

\begin{proposition}
 Let $\mathrm{Arn}$ be the arrow-reversing and Arnold vector subspace of $G$.
 Then 
 $d_{]\cdot[}\bigl(s\mathrm{Arn}\bigr) \subset 
	\bigl(s\mathrm{Arn}\bigr)\cdot\bigl(sG\bigr).$
\end{proposition}

Note that graded anti-commutativity of the graph cobracket in $G$ corresponds via 
$d_{]\cdot[}$ 
to graded commutativity of multiplication in
$\freeA sG$.
$$\xymatrix@R=0pt@C=0pt{
\ \ \ \displaystyle ]g[ & = &
   \displaystyle \sum_e g^{\hat e}_1\otimes g^{\hat e}_2 & = &
   \,\phantom{\frac{1}{2}}
  \displaystyle \sum_e -(-1)^{|g^{\hat e}_1||g^{\hat e}_2|}\, 
                g^{\hat e}_2\otimes g^{\hat e}_1 \qquad\quad\ \\
\displaystyle d_{]\cdot[}\, sg\  & = &
  \displaystyle 
   {\textstyle \frac{1}{2}}\sum_e (-1)^{|g^{\hat e}_1|}\, 
       sg^{\hat e}_1\cdot sg^{\hat e}_2
     & = &
  \displaystyle 
   {\textstyle \frac{1}{2}}\sum_e -(-1)^{|g^{\hat e}_1||g^{\hat e}_2|+|g^{\hat e}_2|}\,
      sg^{\hat e}_2\cdot sg^{\hat e}_1. 
}$$

\begin{corollary}
  $\hat A$ descends to a well-defined map $\A:\dge \to \dga$ by
  $\A([G]) = \hat \A(G)$. 
\end{corollary}

\subsection{Adjointness of $\G$ and $\hat \A$}\label{SS:EAadjoint}

Let $G$ be a $\dgg$ and $A$ be a $\dga$, 
and use $[-]_\g$ to denote the forgetful functor 
to underlying graded vector spaces 
and $[-]_\gg$, $[-]_\ga$ to denote forgetting only differentials. 
It follows from the adjointness properties of $\cofreeG$ and $\freeA$
that the following spaces of homomorphisms are isomorphic:
\begin{equation}\label{adjointness}
{\rm Hom}_\gg\! \left([G]_\gg,\, \cofreeG s^\inv [\bar A]_\g\right)
    \:\cong\: {\rm Hom}_\g\! \left( [G]_\g, \,  s^\inv [\bar A]_\g\right) 
    \:\cong\: {\rm Hom}_\g\! \left( s[G]_\g,\, [\bar A]_\g\right)
    \:\cong\: {\rm Hom}_\ga\! \left(  \freeA s[G]_\g,\, [A]_\ga \right).
\end{equation}
This establishes adjointness of $\G$ and $\hat \A$ on the 
level of graded commutative algebras
and graded graph coalgebras, forgetting differentials.

To display an adjointness which respects $d_\mu$ and
$d_{]\cdot[}$, we translate the classical argument showing adjointness of bar and 
cobar constructions using twisting functions.  
We include the proof only to underline 
that the classical proof translates perfectly to this setting without any 
modification, even though we are now working with the much larger category
of graph coalgebras. 

\begin{theorem}\label{T:adj}
The functors $\G$ and $\hat \A$ are an adjoint pair.
\end{theorem}

\begin{proof}
Given $G$ and $A$, a $\dgg$ and a $\dga$, we will
say that a degree $-1$ map $\tau:[G]_\g \to [A]_\g$ 
is a {\em twisting function} if it satisfies the requirement
$$d_{\!A}\,\tau + \tau\,d_{\!E} - 
    {\textstyle \frac{1}{2}} \Bigl( \mu_A 
     \circ \bigl(((-1)^{|\cdot|}\tau)\otimes \tau\bigr)\; \circ\;\, 
     ]\!\cdot\![_G\,\Bigr) = 0.$$
We show that there are bijections between $\dga$-maps
$\hat \A G \to A$, $\dgg$-maps $G \to \G A$, and twisting functions
$G \to A$.  In terms of 
Equation~\ref{adjointness}, we show that if $[f]_\ga \in
 {\rm Hom}_\ga\! \left(\freeA s [G]_\g,\, [A]_\ga\right)$
comes from applying the forgetful functor to a map $f \in
{\rm Hom}_\dga\! \left(\hat \A G,\, A\right)$, then the adjoint
map $\tau\in {\rm Hom}_\g\! \left( s[G]_\g,\, [\bar A]_\g\right)$
 is in fact a twisting function.  Furthermore
any $\tau\in {\rm Hom}_\g\! \left( s[G]_\g,\, [\bar A]_\g\right)$ 
 which is also a twisting function will be adjoint
to a map $f\in {\rm  Hom}_\ga\! \left(\freeA s [G]_\g,\, [A]_\ga\right)$
in the image of the forgetful functor from ${\rm Hom}_\dga\!(\hat \A E,\, A)$.  This
will complete one half of the argument.  
The half of the argument for homomorphisms $G \to \G A$ 
is similar.

Let $f:\hat \A G \to A$ and write $\tau:s[G]_\g \to [\bar A]_\g$ for
the adjoint of $[f]_\ga$.  Note that $\tau = [f]_\ga \circ i$ where 
$i$ is the injection map $i:s[G]_\g \hookrightarrow \freeA s[G]_\g$. 
The requirement that $d_A\, f = f\, d_{\hat \A G}$ ensures that 
$\tau$ gives a twisting function.
Explicitly, let $sg\in s[G]_\g$, then 
\begin{align*}
0 &= d_A\, f\, i(sg) - f\, d_{\hat \A G}\, i(sg) \\ 
 &= d_A\, f\, i(sg) - f\, (d_{\freeA s G} + d_{]\cdot[})\, i (sg) \\
 &= d_A\, \tau(sg) -  f\, (- s\, d_G g) - f\, 
      \Bigl({\textstyle \frac{1}{2}} 
         \sum_e (-1)^{|g^{\hat e}_1|}\, s g^{\hat e}_1 \cdot s g^{\hat e}_2\Bigr) 
 \qquad \text{where} \quad ]g[\  = 
         \sum_e g^{\hat e}_1 \otimes g^{\hat e}_2\\
 &= d_A\, \tau(sg) +  \tau(s\, d_G g) -  
      {\textstyle \frac{1}{2}} 
        \sum_e (-1)^{|a_i|}\, \tau(s g^{\hat e}_1)\cdot \tau(s g^{\hat e}_2).
\end{align*}

Conversely, let $\tau:s[G]_\g \to [\bar A]_\g$ give a twisting function
$G\to A$ and let $f: \freeA s[G]_\g \to A$ be the adjoint of $\tau$ given by free 
extension.  To show that $d_A\, f = f\,(d_{\freeA sG} + d_{]\cdot[})$ 
it is enough to check on generators
$sg \in \freeA s[G]_\g$.  On generators we have
\begin{align*}
 d_A\, f (sg) &= d_A\, \tau (sg) \\
 f\, d_{\freeA sG} (sg) &= - \tau( s\, d_{G} g) \\ 
 f\, d_{]\cdot[} (sg) & = f\Bigl({\textstyle \frac{1}{2}}
   \sum_e (-1)^{|g^{\hat e}_1|} s g^{\hat e}_1 \cdot s g^{\hat e}_2\Bigr), 
 \qquad \text{where} \quad ]g[\  = \sum_e g^{\hat e}_1 \otimes g^{\hat e}_2\\
  &= {\textstyle \frac{1}{2}} 
   \sum_e (-1)^{|g^{\hat e}_1|}\, \tau(s g^{\hat e}_1) \cdot \tau(s g^{\hat e}_2).
\end{align*}
However, since $\tau$ is a twisting function we know that
$$d_A\, \tau(sg) + \tau (s\, d_G g) - {\textstyle \frac{1}{2}} 
  \sum_e (-1)^{|g^{\hat e}_1|}\, \tau(sg^{\hat e}_1) \cdot \tau(sg^{\hat e}_2) = 0.$$ 
Substitution yields the desired equality.

We only sketch the  bijection between $G\to \G A$ and 
twisting functions $G \to A$, since it is  given similarly.
Let $f:G\to \G A$ and write 
$\tau:[G]_\g \to s^\inv[\bar A]_\g$ for the adjoint of 
$[f]_\gg: [G]_\gg \to \cofreeG s^\inv [\bar A]_\g$.
Note that $\tau = \pi \circ [f]_\gg$ where $\pi$ is the projection map 
$\pi: \cofreeG s^\inv [\bar A] \twoheadrightarrow s^\inv [\bar A]$.  
By direct computation, the requirement that $\pi\, f\, d_G = \pi\, d_{\G A}\, f$
is equivalent to the condition that $\tau$ is a twisting function.
\end{proof}

The adjointness of our duals of Quillen's functors $\eL$ and $\C$ now follows.

\begin{corollary}
 The functors $\E$ and $\A$ are an adjoint pair.
\end{corollary}

Finally, we summarize our results as follows.

\begin{theorem}
The functor $\E : \dga \to \dge$ factors through the category of differential graded
anti-commutative graph coalgebras.
\end{theorem}

\subsection{Pairings of Quillen functors}

Our graphical approach to the Lie coalgebraic bar construction not only gives
rise to the factorization of the previous section, but allows us to explicitly understand
canonical linear dualities of Lie algebraic and coalgebraic Quillen functors.

\begin{theorem}\label{T:diagram}
The diagram
\begin{equation}\label{bigdi}
\begin{aligned}
\xymatrix@=20pt{
\dgc \ar@<.5ex>[r]^{\eL} \ar@<-.5ex>[d]_{\ast} &
  \dgl \ar@<.5ex>[l]^{\C} \ar@<.5ex>[d]^{\ast} \\
\dga \ar@<-.5ex>[r]_{\E} \ar@<-.5ex>[u]_{\ast} &
  \dge \ar@<-.5ex>[l]_{\A} \ar@<.5ex>[u]^{\ast}
}
\end{aligned}
\end{equation}
displays a duality of adjoint pairs of functors.  In particular, the square sub-diagrams
obtained by starting at any corner and mapping to the opposite  are commutative
up to canonical isomorphism.  
In particular, if $C$ is a differential graded-cocommutative coalgebra which
is linearly dual to a differential graded-commutative algebra $A$, then 
$\E(A)$ is linearly dual to $\eL(C)$ through the
configuration pairing.
\end{theorem}

This result refines the work of Schlessinger-Stasheff by identifying the configuration
pairing as giving rise to the canonical duality between the Lie algebraic and coalgebraic
bar constructions.  

\begin{proof}
We treat separately the commutativity of the squares
which constitute  the theorem.  The first two are restated as follows.

{\em  If $L$ is a differential graded  Lie algebra which
is linearly dual to a differential graded Lie coalgebra $E$, then 
$\C(L)$ is linearly dual to $\A(E)$.}

Write $L = (L_\bullet,\, d_L,\, [\cdot,\cdot])$ and $E = (E^\bullet,\, d_E,\, \;]\cdot[\,)$. 
By definition, we need to establish the duality of the bicomplexes
 $\C(L) = \bigl(\cofreeC s L,\ d_{\cofreeC s L},\ d_{[\cdot,\cdot]}\bigr)$ and
$ \A(E) = \bigl(\freeA s E,\ d_{\freeA s E},\ d_{]\cdot[}\bigr).$
Using standard multiplication/comultiplication duality 
the duality between $L_\bullet$ and $E^\bullet$ induces an algebra/coalgebra 
duality between $\cofreeC sL$ and $\freeA sE$.
Furthermore, since $d_L$ and $d_E$ are linearly dual, their cofree/free extensions
$d_{\cofreeC sL}$ and $d_{\freeA sE}$ will be as well.  It remains to show that 
the maps $d_{[\cdot,\cdot]}$ and $d_{]\cdot[}$ are dual.  However, these are also
cofree/free extensions, namely of the maps 
\begin{align*}
 \cofreeC sL \longrightarrow sL \quad &\text{by}\quad 
    s a \cdot s b \longmapsto (-1)^{|a|} s[a,\, b] \\
 sE \longrightarrow\freeA sE \quad &\text{by}\quad
    s \gamma \longmapsto {\textstyle \frac{1}{2}}\, 
          \sum_e (-1)^{|\gamma^{\hat e}_1|}\, 
               s \gamma^{\hat e}_1\cdot s \gamma^{\hat e}_2
 \qquad \text{where}\quad ]\gamma[\ = \sum_e \gamma^{\hat e}_1\otimes \gamma^{\hat e}_2.
\end{align*}
We verify the duality of these restrictions explicitly, using compatibility
of pairings with our assorted multiplications and comultiplications.  
\begin{align*}
\bigl\langle s\gamma,\ (-1)^{|a|} s[a,\, b] \bigr\rangle 
 &= \bigl\langle \gamma,\ (-1)^{|a|} [a,\, b] \bigr\rangle \\
 &= \bigl\langle ]\gamma[\,,\ (-1)^{|a|} a \otimes b \bigr\rangle  \\
 &= (-1)^{|a|}\sum_e \bigl\langle \gamma^{\hat e}_1,\, a \bigr\rangle\, 
                     \bigl\langle \gamma^{\hat e}_2,\, b \bigr\rangle \\
\Bigl\langle {\textstyle \frac{1}{2}} \sum_e (-1)^{|\gamma^{\hat e}_1|}\,
                    s\gamma^{\hat e}_1\cdot s\gamma^{\hat e}_2,\ 
                   s a\cdot s b \Bigr\rangle 
 &= {\textstyle \frac{1}{2}} \sum_e (-1)^{|\gamma^{\hat e}_1|}
              \bigl\langle s\gamma^{\hat e}_1\otimes s\gamma^{\hat e}_2,\ 
                   \Delta (s a\cdot s b)\bigr\rangle \\
 &= {\textstyle \frac{1}{2}} \sum_e (-1)^{|\gamma^{\hat e}_1|}\Bigl(  
                   \langle s\gamma^{\hat e}_1,\ s a \rangle\, 
                   \langle s\gamma^{\hat e}_2,\ s b \rangle \\   
 &\phantom{= {\textstyle \frac{1}{2}} \sum_i (-1)^{|\alpha_i|}\bigl( \ \ }
      + (-1)^{(|a|+1)(|b|+1)}
                   \langle s\gamma^{\hat e}_1,\ s b \rangle\, 
                   \langle s\gamma^{\hat e}_2,\ s a \rangle \Bigr) \\
 &= \sum_e (-1)^{|\gamma^{\hat e}_1|}\langle 
    s\gamma^{\hat e}_1,\ s a \rangle\, \langle s\gamma^{\hat e}_2,\ s b \rangle
\end{align*}
The equality of the last two lines above uses anti-cocommutativity of the cobracket $]\gamma[$ as 
well as the fact that, for the pairings to be nonzero, the degrees of  
$\gamma^{\hat e}_1$ and $a$ must match, as must the degrees of $\gamma^{\hat e}_2$ and $b$.

Since each of the above pairings are $0$ unless $|\gamma^{\hat e}_1| = |a|$, we have equality,
establishing the first half of the theorem.  

The proof of the second half of the theorem
proceeds in the same manner as that of the first half.
Briefly, if we write $A = (A^\bullet, \, d_A, \, \mu)$ and $C = (C_\bullet,\, d_C,\, \Delta)$,
then the duality of the bicomplexes defining $\E(A)$ and $\eL(C)$
is immediate, given by the configuration pairing as stated, 
except for that of the differentials $d_\mu$ and $d_\Delta$.
But  $d_\mu$ and $d_\Delta$ are also cofree/free extensions, namely of the maps
\begin{align*}
 \cofreeE s^\inv \bar A \longrightarrow s^\inv \bar A \quad &\text{by}\quad
    \linep{1}{2} {\textstyle \bigotimes}\, s^\inv a \otimes s^\inv b  
              \longmapsto (-1)^{|a|}\, s^\inv(ab) \\
 s^\inv \bar C \longrightarrow \freeL s^\inv \bar C \quad &\text{by}\quad
    s^\inv \gamma \longmapsto \sum_i (-1)^{|\alpha_i|}[s^\inv \alpha_i,\, s^\inv \beta_i] 
   \qquad \text{where}\quad \bar \Delta \gamma = \sum_i \alpha_i \otimes \beta_i 
\end{align*}
The  duality of these restrictions follows from direct calculation, as before.
\end{proof}

Note that the statements given in the previous proof
do not require our underlying finiteness hypotheses.
If we start with a linearly dual pair of an algebra and coalgebra, the functors $\eL$
and $\E$ will produce a linearly dual Lie algebra and coalgebra.  The finite generation
hypotheses only ensure that our vertical linear duality maps are isomorphisms.

\appendix

\section{Application to computing rational homotopy groups}

We will now collect a number of facts and constructions that were either in the
literature (Schlessinger-Stasheff, Bausfield-Gugenheim) or were  ``in the air''  during
the formative years of rational homotopy theory.  
We are starting to see that a significant pay-off will be obtained when moving to the non-simply-connected case, where our graph coalgebra approach can give rise to additional understanding of 
fundamental groups themselves, rather than having the fundamental group act
on a (minimal) model.   Such results will be the focus of future work.  For the sake of
reference, we collect  first results in the simply-connected setting here.

We discovered the functor $\E$ in the process of  
defining  functionals on homotopy groups,
which in the literature are referred to as homotopy periods.  
Combining our results with the standard translation from spaces to
differential graded algebras shows that this formalism 
is a perfect setting for homotopy periods.
Let $A_*(X)$ be a $\dgc$ model for the rational space $X$, most often
given by the $PL$ chains functor \cite{FHT01}. 
By Quillen's theorem, we know that $H_*(\eL(A_*(X)))$
is isomorphic to $\pi_*(X) \otimes \Q$.  Let $A^*(X)$ be the linear dual
to $A_*(X)$, in other words the $PL$ cochains functor, and let
$\pi_\mathbb{Q}^*(X) = {\rm Hom}(\pi_*(X), \Q)$.

\begin{corollary} \label{C:dualpi}
The homology of $\E(A^*(X))$ is isomorphic to $\pi^*_\mathbb{Q}(X)$.
\end{corollary}

The standard way to recover homotopy data from cochains to this point
has been essentially to replace $A^*(X)$ with a quasi-isomorphic $\A(E)$ 
for some Lie coalgebra $E$, from which it follows by Quillen's theorem
that $\pi^*(X) \cong E$ (see also \refC{mmgood} below).  
Our approach has a number of properties which will
be useful in some settings.

 In the sequel to this paper \cite{SiWa07}, we develop geometry
underlying \refC{dualpi}, defining homotopy periods for any cycles 
in $\E(A_{PL}^*(X))$.  This geometry unifies and generalizes approaches
of Hopf, Whitehead, Boardman-Steer, Sullivan, Novikov, Chen and Hain, 
and can yield $\Z$ and $\Z/p$-valued
homotopy periods.  

Finally, we may employ the spectral sequence of a bicomplex, which  yields the following.

\begin{corollary}\label{C:pi*SS}
If $X$ is a finite complex, 
there is a spectral sequence converging to $\pi^*(X)$ with $E^1$ given by
$\E(H^*(X))$.  This spectral sequence collapses at $E^2$ if $X$ is formal.
\end{corollary}

After \refC{mmgood} we show that this spectral sequence is 
isomorphic to one constructed by Halperin and 
Stasheff \cite{HaSt79} using deformations of minimal models.

\section{Model structures}

We now note that the adjointness results of Section~\ref{SS:EAadjoint}
preserve model structures, so that $\E$, $\A$ and also 
$\G$, $\hat \A$ form Quillen
adjoint pairs. 
Because we are in the finitely generated setting, we get only 
model structures, not closed model structures.  

All categories in this section are 
reduced appropriately.

\begin{theorem}[Quillen \cite{Quil69}]\label{T:LCmodel}
 A model category structure on $\dgl$ is given by the following: 
\begin{itemize}
\item Weak equivalences are the quasi-isomorphisms. 
\item Fibrations are the level-wise surjections above the bottom degree.
\item Cofibrations are determined by left lifting; they are 
the free $\gl$-maps.
\end{itemize}

 A  model category structure on $\dgc$ is given by the following: 
\begin{itemize}
\item Weak equivalences are the quasi-isomorphisms. 
\item Cofibrations are the levelwise injections.
\item Fibrations are determined by right lifting. 
\end{itemize}
\end{theorem}

Recall that $\otimes$ gives finite products in $\gc$, since our coalgebras are counital, coaugmented.
Note that all $\dgc$'s are cofibrant and all $\dgl$'s are fibrant.

\begin{remark}
By the results of Quillen~\cite{Quil69}, these give model category 
structures even with the finiteness assumptions
removed.  Though Quillen did not show that these model categories are 
closed when finiteness hypotheses are removed,
in particular that infinite limits exist in the coalgebra 
setting, there are now a number of proofs in the literature.
\end{remark}

In the course of developing algebraic models for rational homotopy theory,
Quillen established the following (see \cite[Thm 5.3]{Quil69}).

\begin{theorem}[Quillen]\label{T:QuiAdj}
The functors $\eL:\dgc \leftrightarrows \dgl: \C$ are a Quillen adjoint pair.  That is, 
$\eL$ preserves cofibrations and trivial cofibrations; $\C$ preserves fibrations and
trivial fibrations.

Furthermore, $\eL$ and $\C$ give a Quillen equivalence.  That is, if $C$ is a cofibrant
$\dgc$ and $L$ is a fibrant $\dgl$, then a map $\eL (C) \to L$ is a weak equivalence if and only
if the adjoint map $C \to \C (L)$ is a weak equivalence.
\end{theorem}

We now give parallels to these results in our algebra--Lie coalgebra setting.
In the following, we continue to restrict to finitely generated, reduced objects.

\begin{definition}

We will say that a $\dga$-map $f:A \to B$ is a free $\ga$-map if as a
$\ga$-map, it is an inclusion of a graded algebra 
with free cokernel, as displayed in the diagram:
$$\xymatrix@R=7pt@C=25pt{
[A]_\gl \ar[r]^{[f]_\gl} \ar@{_(->}[dr] 
   &  [B]_\gl \\
   & [A]_\gl \otimes \freeA W .
    \ar@{{}{}{}}[u]^*[right]{\cong}
}$$
In $\ga$, $\otimes$ is the categorical coproduct, since our algebras are unital.

We will say that a $\dgg$-map $f:D\to E$ is a cofree $\gg$-map if as a  
$\gg$-map, it is a projection of graded coalgebras with 
cofree kernel, as displayed in the diagram:
$$\xymatrix@R=7pt@C=25pt{
[D]_\gg \ar[r]^{[f]_\gg} 
  \ar@{{}{}{}}[d]_*[right]{\cong}
   &  [E]_\gg \\
[E]_\gg \circledast \cofreeG W. \ar@{->>}[ur] &
}$$
By  $\circledast$ we mean the ``cofree product'' -- the categorical product of 
graph coalgebras --
given by the categorical equalizer of the pair of maps
$$G\circledast K := \mathrm{Eq}
 \Bigl(\cofreeG(G\oplus K) \rightrightarrows \cofreeG(\cofreeG G \oplus \cofreeG K)\Bigr)$$
coming from $\cofreeG$ being a cotriple and from $G$, $K$ being graph coalgebras.

\end{definition}

\begin{theorem}\label{T:AEmodel}
 A model category structure on $\dga$ is given by the following:
\begin{itemize}
\item Weak equivalences are the quasi-isomorphisms.
\item Fibrations are the levelwise surjections.
\item Cofibrations are determined by left lifting; they are the free $\ga$ maps.
\end{itemize} 

 A model category structure on $\dgg$ is given by the following:
\begin{itemize}
\item Weak equivalences are the quasi-isomorphisms.
\item Cofibrations are injections above degree one.
\item Fibrations are determined by right lifting; they are the cofree $\gg$ maps.

 A model category structure on $\dge$ is given similarly.
\end{itemize}

\end{theorem}

While it is possible to merely mimic the original proof of Quillen from~\cite{Quil69},
we may instead infer this from the literature on model categories.

\begin{proof}[Proof Sketch]

The stated model category structure on $\dga$ is standard in the literature -- it is
given by lifting the projective model structure on (reduced) cochains.  See \cite{Hirs03}
and \cite[4.1]{ScSh00}.

To see that the cofibrations are indeed the free maps may be done in the same way as
Quillen shows the corresponding fact in $\dgl$ (see~\cite[Prop 5.5, p256]{Quil69})
by attatching cells using pushouts of cofibrations.  In this manner one may show that
all cofibrations are retracts of free maps.  However, subalgebras
of free algebras are again free; so such maps must themselves be free.  

The listed model category structure on $\dgg$ 
is implied by general operad theory work of~\cite[Thm 3.2.3]{AuCh03}.
That fibrations are indeed the cofree maps follows in the finitely generated case from
the dual of the corresponding statement about cofibrations in $\dgl$. 
\end{proof}

\begin{remark}
As in the $\dgc$ and $\dgl$ settings, 
the structures given in Theorem~\ref{T:AEmodel} (minus the 
description of fibrations in $\dgg$) give closed model category structures when
finiteness assumptions are removed. 
There is a discrepancy between this situation and that of 
\cite{AuCh03}, which defines
cooperads using direct sums and orbits instead of products and fixed points.
\end{remark}

\begin{lemma}\label{L:dualmodel}
The model structures of $\dgl$ and $\dge$ and of $\dga$ and $\dgc$ given 
in \refT{LCmodel} and \refT{AEmodel} are linearly dual.  
That is, each vertical linear
duality isomorphism sends fibrations to cofibrations, cofibrations to fibrations,
and weak equivalences to weak equivalences.
\end{lemma}

While we have generally chosen to give self-contained arguments, for 
showing that $\E$ and $\A$ give a Quillen equivalence
we stray from this choice for the sake of brevity.  
We may deduce the following result from 
\refL{dualmodel}, our main \refT{diagram}, 
and Quillen's Theorem as stated in \refT{QuiAdj}.

\begin{theorem}\label{T:EAQuiAdj}
The functors $\G$ and $\hat \A$ are a Quillen adjoint pair.

The functors $\E$ and $\A$ are a Quillen adjoint pair.  Further, $\E$ and $\A$ are 
a Quillen equivalence.
\end{theorem}

\section{Minimal models}\label{S:minimal models}

We end with some brief notes about minimal models, originally due to Sullivan
 \cite{Sull77, FHT01}.   In our language, 
a minimal model in $\dga$ is an object of the form $(\freeA W,\,d)$ where
$d\,W\subset\freeA^{\ge 2} W$.  
Sullivan's theorem \cite{Sull77} is that every 
$\dga$ supports a quasi-isomorphism from a minimal model
$(\freeA W,\ d) \xrightarrow{\ \simeq\ } A$, and furthermore the minimal
model $(\freeA W,\ d)$ is unique up to isomorphism. 
Minimal models in $\dga$ are
useful because the Postnikov tower of a rational one-reduced
space is encoded transparently in its minimal model as the 
increasing filtration by free sub-algebras. 

Baues and Lemaire \cite{BaLe77} note that the property satisfied by the differential of
a minimal model may be more concisely stated as $(-)^{\text{ind}} \circ d = 0$.  Further,
they show that making the analogous definition in $\dgl$ also agrees with the naive definition,
namely $(\freeL V,\, d)$ with $d\,V\subset\freeL^{\ge 2}V$.
These minimal models have existence and uniqueness properties similar to 
those of Sullivan's minimal models in $\dga$, but because of the switch
 from cochains to chains their construction is 
more difficult -- see \cite{BaLe77}.  From the point of view of topology,
minimal models in $\dgl$ encode the 
Eckmann-Hilton homology decomposition of a rational space.

One lemma in the proof of the uniqueness of minimal models of algebras 
is interesting in its own right.  We say that a $\dga$ is a ``differential 
free graded algebra'' if it has the form $(\freeA V,\, d)$, and similarly for a
``differential free graded Lie algebra''.  Then we have the following \cite[Prop 1.5]{BaLe77}.

\begin{proposition}[Sullivan, Baues-Lemaire]\label{P:Q-Hurewicz}
A map $f$ of differential free graded (Lie) algebras is a quasi-isomorphism 
if and only if the induced $\dg$-map $(f)^{\text{ind}}$ on indecomposables is
a quasi-isomorphism.
\end{proposition} 

We apply this proposition to the units of the adjunctions $\A\, \E \to \mathbbm{1}_\dga$ and 
$\eL\, \C \to \mathbbm{1}_\dgl$.

\begin{corollary}\label{C:mmgood}
If $A$ is a differential free graded algebra, then $[\E A]_\dg \simeq s(A)^{\text{ind}}$.
Similarly, if   $L$ is a 
differential free graded Lie algebra, then  $[\C L]_\dg \simeq s(L)^{\text{ind}}$.

In particular if $A$ is a $\dga$ minimal model, then 
$H^* \E A \cong s(A)^{\text{ind}}$ as a graded vector space.
Similarly, if $L$ is a $\dgl$ minimal model then $H_* \C L \cong s(L)^{\text{ind}}$.
\end{corollary}

We use this corollary to recover the Halperin-Stasheff 
spectral sequence for calculating the linear dual of homotopy groups of a finite
complex, as described in 4.14 of \cite{HaSt79}, from our \refC{pi*SS}.  
The main construction of \cite{HaSt79} is that of a filtered model
for $(A,\,d_A)$ as a deformation of a minimal model for 
$(H_{*}(A),\, 0)$, which in our notation would be 
called $(\freeA Z,\, D)$ and $(\freeA Z,\, d)$ respectively.    When $A = A^{*}(X)$,
the results of Section 8 of \cite{Sull77} imply that $H_{*}(Z,\, D) \cong \pi^{*}(X)$.  Because
$D$ and $d$ differ by terms of lower filtration, there is a spectral sequence 
starting with $H_{*}(Z,\, d)$ and converging to $H_{*}(Z,\, D) \cong \pi^{*}(X)$.  

By \refC{mmgood}, we have $H_{*}(Z,\, d) \cong H_{*}\left(\E(H^{*}(X)\right)$,
so this spectral sequence has the same $E^{2}$ term as that of \refC{pi*SS}.
Indeed,  we may relate these two spectral sequences by comparing them both
to equivalent spectral sequences for $\E(\freeA Z,\, D)$, which on one hand is quasi-isomorphic to 
$\E(A,\, d_A)$ simply because $\E$ is quasi-isomorphism invariant; and on the other hand is
quasi-isomorphic to $(Z,\, D)$ by \refC{mmgood}.  Our approach through $\E(A,\, d_A)$
seems to have  better functorality properties, a more transparent cobracket structure,
and greater flexibility in addition to the conjectured relationship with Hopf invariants.

\medskip

Natural notions of minimal models in coalgebras are obtained by duality.  Explicitly
we require them to be cofree with differentials satisfying $d\circ (-)^{\text{pr}} = 0$.

\begin{definition}
A minimal model in $\dgc$ is a coalgebra of the form $(\cofreeC V,\ d)$ where
$d\,V=0$. 

A minimal model in $\dge$ is a coalgebra of the form $(\cofreeE W,\ d)$ where
$d\,W=0$.
\end{definition}

We may speak of ``differential cofree graded (Lie) coalgebras'' similarly to obtain
duals to Proposition~\ref{P:Q-Hurewicz} and Corollary~\ref{C:mmgood}.

\begin{proposition}
A map $f$ of differential cofree graded (Lie) coalgebras is a quasi-isomorphism
if and only if the induced $\dg$-map $(f)^{\text{pr}}$ on primitives is a 
quasi-isomorphism.
\end{proposition}

\begin{corollary}\label{C:indec}
If $C$ is a differential cofree graded coalgebra, then $[\eL C]_\dg \simeq s^\inv(C)^{\text{pr}}$.
Similarly, if $E$ is a differential
cofree graded $\eil$ coalgebra, then $[\A E]_\dg \simeq s^\inv(E)^{\text{pr}}$.

In particular if $C$  is a $\dgc$ minimal model, then 
$H_* \eL C \cong s^\inv (C)^{\text{pr}}$ 
as a graded vector space.  Similarly, if $E$ is a $\dge$ minimal model,
then  $H^* \A E \cong s^\inv (E)^{\text{pr}}$.
\end{corollary}

Minimal models in all cases are unique up to isomorphism for each object, 
an Bousfield-Gugenheim even give a functorial construction of them \cite{BoGu76}.
Minimal models of algebras are cofibrant replacements, and minimal models of coalgebras
are fibrant replacements.  There are other standard functorial fibrant and cofibrant 
replacements, namely in each setting by applying the appropriate pair of adjoint horizontal 
arrows from the diagram of \refT{diagram}.
These generally differ from minimal models, and as indicated by our discussion of the 
Halperin-Stasheff spectral sequence the interplay between the two approaches can be
enlightening.


\end{document}